\documentclass[twocolumn]{autart}    %

\usepackage{hyperref}
\usepackage{amsmath}
\usepackage{amsfonts}
\usepackage{mathtools}
\usepackage{bm}
\usepackage{enumitem}
\usepackage{accents}
\usepackage{float}
\usepackage{graphicx}
\usepackage{makeidx}
\usepackage[caption=false,font=footnotesize]{subfig}
\usepackage[cal=cm,calscaled=1.0,bb=ams,frak=pxtx,frakscaled=1.0,scr=boondox,scrscaled=1.0,bb=boondox]{mathalfa}

\usepackage{xcolor}
\usepackage{amssymb}
\usepackage{nicefrac}

\usepackage{breqn}

\usepackage{physics}

\usepackage{natbib}

\newcommand\Paverage[1]{\mathbb{P}\left(#1\right)}

\newcommand\condExpectA[2]{\mathbb{E}\left[\left.#1\right|#2\right]}
\newcommand\rhoConditional[2]{\mathbb{E}_{#1}\left[#2\right]}

\usepackage{etoolbox}

\makeatletter
\preto\maketitle{%
	\begingroup\lccode`~=`,
	\lowercase{\endgroup
		\let\saved@breqn@active@comma~%
		\let~}\active@comma %
}
\appto\maketitle{%
	\begingroup\lccode`~=`,
	\lowercase{\endgroup
		\let~}\saved@breqn@active@comma %
}
\makeatother
\usepackage[toc]{glossaries}
\newacronym{QSDE}{QSDE}{quantum stochastic differential equation}
\newacronym{MBFC}{MBFC}{measurement-based feedback control}

\newacronym{KF}{KF}{{Kalman} filter}
\newacronym{EKF}{EKF}{extended {Kalman} filter}
\newacronym{quantum EKF}{quantum EKF}{quantum extended {Kalman} filter}
\newacronym{MIMO}{MIMO}{multiple-input multiple-output}
\newacronym{SME}{SME}{stochastic master equation}
\newacronym{QDS}{QDS}{quantum dynamical semigroup}
\newacronym{QMS}{QMS}{quantum {Markovian} semigroup}
\newacronym{SSE}{SSE}{stochastic {Schr\"{o}dinger} equation}
\newacronym{SDE}{SDE}{stochastic differential equation}
\newacronym{HJI}{HJI}{{Hamilton} {Jacobi} inequality}
\newacronym{MIST}{MIST}{monotonically increasing sequence of times}
\newacronym{LMI}{LMI}{linear matrix inequality}
\newacronym{MISE}{MISE}{mean integral of squared quadratic error}
\newacronym{OPD}{OPD}{ordered partial differentiation}
\newacronym{OD}{OD}{ordered differential}
\newacronym{OM}{OM}{ordered multiplication}

\begin{document}

\begin{frontmatter}

\title{Stability Analysis of Quantum Systems: a Lyapunov Criterion and an Invariance Principle} %

\thanks[footnoteinfo]{This work was supported by the Australian Research Council under grant numbers FL110100020 and DP180101805 and the AFOSR under agreement number FA2386-16-1-4065.}

\author[ADFA,ANU,AALTO]{Muhammad F. Emzir}\ead{muhammad.emzir@aalto.fi},    %
\author[ADFA,PRINCETON]{Matthew J. Woolley}\ead{m.woolley@adfa.edu.au},               %
\author[ANU]{Ian R. Petersen}\ead{i.r.petersen@gmail.com}  %
\address[ADFA]{School of Engineering and Information Technology, University of New South Wales, ADFA, Canberra, ACT 2600, Australia}
\address[PRINCETON]{Department of Electrical Engineering, Princeton University, Princeton, New Jersey 08544, USA}
\address[ANU]{College of Engineering and Computer Science, Australian National University, Canberra ACT 0200, Australia}
\address[AALTO]{ Department of Electrical Engineering and Automation, Aalto University, Espoo 02150, Finland}


\begin{abstract}                          %
In this article, we propose a Lyapunov stability approach to analyze the convergence of the density operator of a quantum system. In analog to the classical probability measure for Markovian processes, we show that the set of invariant density operators is both closed and convex. We then show how to analyze the stability of this set via a candidate Lyapunov operator. We complete our analysis of the set of invariant density operators by introducing an analog of the Barbashin-Krasovskii-La Salle theorem on the dynamics of quantum systems.  %
\end{abstract}

\end{frontmatter}

\section{Introduction}
Among the most fundamental issues in dynamics and control theory is the stability of systems. The Lyapunov criterion is one of the key tools used to analyze stability. Since its discovery in the late nineteenth century\citep{Lyapunov1892}, the Lyapunov criterion has been successfully applied in nonlinear stability analysis and feedback control design. The criterion has also been extended in various ways. Among these extensions are the invariance principle, due to Barbashin-Krasovskii-La Salle \citep{Salle1960,Salle1961,Krasovsky1963}, where the Lyapunov criterion has been extended to infer asymptotic stability when the time derivative of the Lyapunov function fails to be strictly negative in the neighborhood of an equilibrium point. There are also other extensions of the Lyapunov criterion that cover the stability analysis of systems with time delay \citep{Kharitonov1999,Chen2003,Zhou2016,Liu2016,Mazenc2017}, switched systems \citep{Lin2009,Jungers2017}, and stochastic systems \citep{Kozin1969,khasminskii2011stochastic} to mention a few examples. 
\\
The development of feedback control for the quantum systems has received a significant amount of interest in the last two decades \citep{Dong2010}.  Essentially, there are two approaches to feedback control design for quantum systems. The first approach relies on measurements of the system, where the feedback control is calculated based on the estimated value of the variables of the system. This method, known as measurement-based feedback control, has been well investigated in the last two decades \citep{Wiseman1993,Belavkin1999,wiseman2010quantum}. The other approach is to construct the feedback controller as a quantum system which coherently interacts with the controlled system, also known as coherent feedback control. This approach has recently received considerable attention \citep{Wiseman1994,James2008,Nurdin2009,Maalouf2011,Zhang2012}. In some recent studies \citep{Nurdin2009,Hamerly2012,Yamamoto2014}, it has been shown that there are some conditions in which coherent feedback control can offer improvements over measurement-based feedback control .
\\
Regarding analytical tools to analyze the convergence of quantum system dynamics, there have been many results on analytical tools to analyze the convergence of quantum system dynamics based on the stability analysis of quantum systems subject to measurement \citep{somaraju2013,Amini2013,amini2014,Benoist2017}. However, in the absence of measurement, as in the case of coherent control, there are few tools available to analyze the stability behavior of quantum dynamical systems. Many established results on linear coherent control of quantum systems are based on stochastic stability criteria involving the first and second moments \citep{James2008}. To extend coherent control design beyond the linear case, it is useful to consider more general stability criteria other than first and second moment convergence.
\\
From the classical probability theory point of view, we can consider a density operator of the system Hilbert space in a similar way to a probability measure. Therefore, the convergence of a density operator can be analyzed in a similar way to the convergence of a sequence of probability measures. The convergence of probability measures in metric spaces has been studied since the mid-1950's, starting from the seminal work of Prokhorov \citep{Prokhorov1956,Parthasarathy1967}. In fact, in the mathematical physics literature, the stability of quantum systems has been analyzed using semigroups and the quantum analog of probability measure convergence \citep{Fagnola2001,Fagnola2003,Fagnola2006,Chang2014,Deschamps2016}.  In essence, \citep{Fagnola2001,Fagnola2002,Fagnola2003,Fagnola2006} establish conditions for the existence of invariant states and the convergence to these states given that the invariant state is faithful; i.e., for any positive operator $A$, $\trace\qty(A \rho) = 0$ if and only if $A=0$. However, although the quantum Markov semi-group class considered in \citep{Fagnola2001,Fagnola2002,Fagnola2003,Fagnola2006} is quite general, their result is only valid for faithful invariant states. Moreover, the convergence conditions are based on verifying whether a set of operators that commute with all coupling operators and their adjoint, $\qty{\mathbb{L}_i,\mathbb{L}_i^\ast, \forall i}'$, coincides  with a set of operators that commute with the Hamiltonian, all coupling operators and their adjoints, $\qty{\mathbb{H}, \mathbb{L}_i,\mathbb{L}_i^\ast, \forall i}'$, which can be difficult to verify; see \citep[Theorem 6.10]{Fagnola2003}\citep{Arnold2008,Nigro2019}.
\\ 
Stability analysis in the Heisenberg picture as given in \citep{Pan2014} is interesting for two reasons. The first is that since it is written in Heisenberg picture, which directly describes the dynamics of system observables, the stability condition derived is easily connected to classical Lyapunov stability analysis. The second is that while the stability condition is stated in term of a Lyapunov observable, it leads to the same conclusion as the quantum Markov semigroup convergence approach. 
The result of \citep{Pan2014} required that for all non-trivial projection operators $P$,  $P\mathbb{L}^\dagger\qty(\mathbb{1}-P)\mathbb{L}P \neq 0$, where $\mathbb{L}$ is any coupling operator of the quantum system.
Moreover, the crucial assumption of this result and those of \citep{Fagnola2001,Fagnola2002,Fagnola2003,Fagnola2006} is that the steady state (or invariant state) of the evolution of the system density operator is faithful.
In many practical applications, this is not the case, such as in the case of ground state stabilization; see \citep{Mamaev2017} for an example of entangled state dissipation, and other pure state stabilization problems \citep{Yamamoto2012,Ticozzi2012,Leghtas2013,Camalet2013}. If the set of invariant density operators are the ground states of a self-adjoint operator, then these density operators will not be faithful unless the self-adjoint operator is a constant multiple of the identity. This shows that there are problems in stability analysis for quantum dynamical semigroups which cannot be solved by the currently available methods. In fact, convergence analysis to the set of ground states has been addressed partially for the finite dimensional case recently in \citep{Pan2016}.
Furthermore, even when the invariant density operators are faithful, validating the inequality as in \citep[Theorem 3, Theorem 4]{Pan2014} for all non-trivial projection operators is not straight forward; see also \citep[Example 4]{Pan2014}.  
In \citep{Pan2016}, the convergence of the system density operators to a \emph{ground states} of a self adjoint operator $X$, $\mathscr{Z}_X$ is considered for finite dimensional quantum systems. It was shown in \citep{Pan2016} that if for a Lyapunov operator $V$, the generator $\mathcal{L}(V)$ of $V$ has a strictly negative value $\trace\qty(\mathcal{L}(V) \rho)$ for any $\rho \notin \mathscr{Z}_X$, then the system density operator will converge to $\mathscr{Z}_X$. The crucial assumptions made in \citep{Pan2016} are that the Lyapunov operator $V$ has to be a positive operator $V$ satisfying $\mathcal{L}(V) \leq 0$ and has to commute with the Hamiltonian of the system. 
\\
A more general approach to stability analysis of the invariant density operators for quantum systems can be established by considering the convergence of density operators in Banach space, and requiring less restrictive assumptions on the Lyapunov operators.
\\
In this article, we aim to establish a stability criterion, which is similar to the Lyapunov stability theory in classical dynamics, to examine the convergence of a density operator of a quantum system. In contrast to many previously studied convergence analysis methods for invariant density operators which used weak convergence \citep{Fagnola2001,Fagnola2002,Fagnola2003,Fagnola2006}, in this work, we will analyze the convergence of density operators in a Banach space using the Lyapunov criterion. This approach may be more useful in some cases than the previous approaches since the weak topology does not give a concrete distance between two density operators. Hence it is impossible to determine the distance of the density operator from the steady state (invariant density operator). On the other hand, analyzing convergence in the Banach space allows a direct connection to classical stability analysis by generalizing the finite-dimensional vector space into a general Banach space.	It is also worth mentioning that the convergence of density operator evolution has also been analyzed via Wigner-Fokker-Planck partial differential equations, where in \citep{ARNOLD2012}, an exponential convergence rate towards the invariant state is established for a class of perturbed linear quantum stochastic differential equations (QSDEs). Notice however, our stability criterion should \emph{not} be regarded as an alternative to the ergodic theory of quantum dynamical semigroups, but rather as a complementary approach that applies to non-faithful invariant density operators.
\\
After giving a stability definition in the set of density operators, first, we will demonstrate that if there is a self-adjoint operator that has a strict minimum value at the invariant state and its generator satisfies a particular inequality, then we can infer Lyapunov, asymptotic, and exponential stability in both local and global settings. We do not require the Lyapunov operator to be either positive or to be coercive (the spectrum of $V$ is greater than some strictly increasing function and goes to infinity), as considered in previous work \citep{Pan2014}. This flexibility enables us to study the stability of quantum systems with non-faithful invariant density operators. 
\\
We complete our analysis of the set of invariant density operators by introducing a quantum analog of the Barbashin-Krasovskii-LaSalle Theorem. In the classical setting, this theorem states that the trajectory of a system will approach the greatest invariant set of the system when there is a Lyapunov function which is non-increasing over time \citep{Krasovsky1963,khalil2002nonlinear}. 
The main departure of the quantum result from the classical Barbashin-Krasovskii-La Salle Theorem is that we require that the solution of the predual semi-group $\rho_t$ belongs to a weak compact set.  Note that the different topologies only matter if we deal with an infinite dimensional Hilbert space, since weak and normal convergence are equivalent in finite-dimensional normed spaces  \citep[Prop 2.5.13]{Megginson1998}.
\\
Part of this paper has been presented at the 56th IEEE Control and Decision Conference \citep{Emzir2017a}. Compared to the conference version, this paper has been significantly expanded and revised. Section \ref{sec:InvariantPrinciple}, which discusses the invariance principle analog for quantum systems, has been added to the paper.  Further, we have also included Theorem \ref{thm:Global_Local_Asymptotic} which shows that local asymptotic stability implies global asymptotic stability in the quantum systems being considered.  In addition to that, Example \ref{sec:exm_LinearQuantumSystemWithOtherV} has been added to illustrate the weak asymptotic stability condition via our quantum version of the Barbashin-Krasovskii-La Salle Theorem. The paper also contains complete proofs of all results, which have not been included in the conference paper \citep{Emzir2017a}.
\\
Following a common approach in the physics literature, we will base our work on the algebras of bounded operators. 
Notice however, most unbounded operators used in this work (e.g. position and momentum operators of a quantum particle) have domains that include at minimum the set of smooth functions with compact support that are dense in $\mathcal{L}_2\qty(\mathbb{R})$ \citep{Parthasarathy1992,bouten2007introduction}.
We refer the reader to \citep[\textsection 2.1]{Emzir2018a} for a brief introduction to the closed and open quantum system dynamics used in this article. We also refer the reader to the following monographs for an introduction to quantum probability \citep{bouten2007introduction,fagnola1999quantum}. 

\subsection{Notation}
The symbol $\mathbb{1}$ stands for the identity operator of a Hilbert space. Plain letters (e.g., $P$) will be used to denote operators on Hilbert space. Script face (e.g., $\mathscr{H}$ for Hilbert space) is used to denote spaces as well as algebras. A class of operators will be denoted by fraktur face; e.g., we let $\mathfrak{B}\left(\mathscr{H}\right)$ to denote the class of bounded linear operators on the Hilbert space $\mathscr{H}$. $\mathfrak{I}_1\left(\mathscr{H}\right)$ denotes the set of trace-class Hilbert space operators. For a bounded operator $X$, the uniform norm is denoted by $\norm{X}_\infty$, while if $X$ is a trace-class operator, the norm $\norm{X}_1$ is defined as $\norm{X}_1 = \trace\qty(\sqrt{X^\ast X})$.
Bold letters (e.g., $\mathbf{y}$) will be used to denote a matrix whose elements are operators on a Hilbert space. Hilbert space adjoints, are indicated by $^{\ast}$, while the adjoint transpose will be denoted by $\dagger$; i.e., $\left(\mathbf{X}^{\ast}\right)^{\top} = \mathbf{X}^{\dagger}$. For single-element operators we will use $*$ and $\dagger$ interchangeably. %

\section{Quantum Dynamical Semigroups}\label{sec:PreliminaryLyapunovStability}
\begin{figure}[h]
	\centering
	\includegraphics[width=0.5\textwidth]{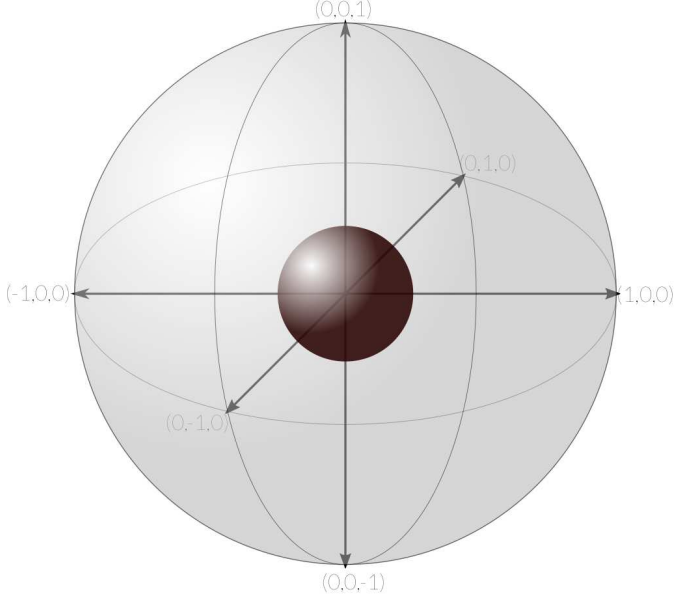}
	\caption{Illustration of a ball in \eqref{eq:BallCenterAtRhoStarRadiusEpsilon} for a qubit system, drawn in the Bloch sphere.} 
	\label{fig:QubitBall}
\end{figure}

In this section, we will describe some notions in quantum dynamical semigroups that will be used in the later sections; see \citep[Chapter 1]{davies1976quantum},\citep[\textsection 2-3]{fagnola1999quantum},\citep[\textsection 2.6]{Chang2015}, \citep[Chapter 1]{Emzir2018}. \\
A class of density operators $\mathfrak{S}\qty(\mathscr{H})$ on the Hilbert space $\mathscr{H}$ is defined as a subset of the trace class operators $\mathfrak{I}_1\qty(\mathscr{H})$ with non-negative spectrum and unity trace, \citep[\textsection 1.4]{davies1976quantum}.
Notice that a linear functional $\omega$ on $\mathfrak{B}\left(\mathscr{H}\right)$ is  normal if there is a unity trace operator $\rho$ such that $\omega(X) = \trace\qty(\rho X)$. From this viewpoint, a unit element $\ket{u} \in \mathscr{H}$ and a density operator $\rho$ can also be considered as states on $\mathfrak{B}\qty(\mathscr{H})$, by considering the following linear functionals, $\varphi\qty(X) =  \bra{\psi}X\ket{\psi}$, and $\varphi\qty(X) = \trace\qty(\rho X)$.
Let $\mathscr{H}$ be the system's Hilbert space. The total Hilbert space will be given as $\tilde{\mathscr{H}} = \mathscr{H} \otimes \Gamma$ where $\Gamma$ is the Fock space for the environment \citep{Emzir2018a}. 
Notice that if we restrict ourself to work in the linear span of coherent states, for any time interval $0\leq s < t $, the Fock space $\Gamma$ can be decomposed into \citep[pp. 179-180]{Parthasarathy1992}
$
	\Gamma = \Gamma_{s]}\otimes \Gamma_{[s,t]}\otimes \Gamma_{[t}. \label{eq:GammaDecompose}
$	
Therefore, we can write $\tilde{\mathscr{H}}_{\left[0,t\right]} = \mathscr{H} \otimes {\Gamma}_{\left[0,t\right]}$.
If $\rho \in \mathfrak{S}\qty(\mathscr{H})$ is an initial density operator of the system and $\Psi \in  \mathfrak{S}\qty(\Gamma)$ is the initial density operator of the environment, then for any bounded system observable $X$, the quantum expectation of $X_t = j_t(X) \equiv U_t^\dagger( X \otimes \mathbb{1})U_t$ is given by $\trace{ \qty(X_t \qty(\rho \otimes \Psi))} \equiv \Paverage{X_t}$. Let $\Psi_{[t} \in  \mathfrak{S}(\tilde{\mathscr{H}}_{[t})$ be a density operator on $\tilde{\mathscr{H}}_{[t}$. Then, we define $\rhoConditional{t]}{\cdot} : \mathfrak{B}\qty(\tilde{\mathscr{H}}) \rightarrow \mathfrak{B}(\tilde{\mathscr{H}}_{t]})$ as follows:
\begin{dmath}
\rhoConditional{t]}{Z} \otimes \mathbb{1} \equiv \condExpectA{Z}{\mathfrak{B}\qty(\tilde{\mathscr{H}}_{t]})\otimes \mathbb{1} } {, \forall Z \in \mathfrak{B}\qty(\tilde{\mathscr{H}})}
\label{eq:PartialTraceExpect}
\end{dmath}
where $\mathbb{1}$ is the identity operator on $\tilde{\mathscr{H}}_{[t}$ and $ \condExpectA{Z}{\mathfrak{B}\qty(\tilde{\mathscr{H}}_{t]})\otimes \mathbb{1} }$ is the quantum conditional expectation; see also \citep[Prop 16.6, Exc 16.10, 16.11]{Parthasarathy1992} and \citep[Exm 1.3]{fagnola1999quantum} for the existence of $\condExpectA{Z}{\mathfrak{B}\qty(\tilde{\mathscr{H}}_{t]})\otimes \mathbb{1} }$. In our case, we will frequently consider $\rhoConditional{0]}{j_t\qty(X)}$, when the quantum expectation of $j_t\qty(X)$ is marginalized with respect to the system Hilbert space $\mathscr{H}$.
\begin{defn}\citep{fagnola1999quantum,Chang2015}
	A quantum dynamical semi-group (QDS) on a von Neumann algebra $\mathscr{A}$ is a family of bounded linear maps $\qty{\mathcal{T}_t, t\geq 0}$ with the following properties:
	\begin{enumerate}
		\item $\mathcal{T}_0\qty(A) = A$, for all $A \in \mathscr{A}$.
		\item $\mathcal{T}_{s+t}\qty(A) =  \mathcal{T}_{s}\qty(\mathcal{T}_{t}\qty(A))$, for all $s,t \geq 0$, and $A \in \mathscr{A}$.
		\item $\mathcal{T}_{t}$ is a completely positive mapping for all $t\geq 0$.
		\item $\mathcal{T}_{t}$ is normally continuous for all $t\geq 0$.
	\end{enumerate} 
\end{defn}
Using the conditional expectation in \eqref{eq:PartialTraceExpect}, we observe that there exists a one-parameter semi-group $\mathcal{T}_t : \mathfrak{B}\qty(\mathscr{H}) \rightarrow \mathfrak{B}\qty(\mathscr{H})$ given by 
\begin{dmath}
{\mathcal{T}_t\qty(X) = \rhoConditional{0]}{j_t\qty(X)} }.\label{eq:TSemiGroup}
\end{dmath}
The generator of this semi-group $\mathcal{L}\qty(X) : \mathscr{D}\qty(\mathcal{L}) \rightarrow \mathfrak{B}\qty(\mathscr{H})$ is given by
$
\mathcal{L}\qty(X) = \lim_{t \downarrow 0} \frac{1}{t}(\mathcal{T}_t\qty(X) - X), \forall X \in \mathscr{D}\qty(\mathcal{L})
$
. By using the quantum conditional expectation property \citep[\textsection 2-3]{fagnola1999quantum} and \eqref{eq:PartialTraceExpect}, we obtain $\Paverage{\mathcal{T}_t\qty(X)\otimes \mathbb{1}} = \Paverage{X_t}$. Therefore,
\begin{dmath}
{\mathcal{T}_t\qty(\mathcal{L}\qty(X)) = \mathcal{L}\qty(\mathcal{T}_t\qty(X)) = \rhoConditional{0]}{\mathcal{G}\qty(X_t)}}; \label{eq:LandGrelation}
\end{dmath}
see also \citep[Thm 5.1.6]{Chang2015}.
Note that there also exists a one-parameter semi-group, known as the predual semi-group $\mathcal{S}_t$ such that, $\Paverage{X_t} = \Paverage{\mathcal{T}_t\qty(X) \otimes \mathbb{1}} = \trace\qty(X \mathcal{S}_t\qty(\rho))$. Explicitly, it can be defined as
\begin{dmath}
\mathcal{S}_t\qty(\rho) \equiv \trace_{\Gamma}\qty(U_t(\rho \otimes \Psi)U_t^\dagger), \label{eq:SSemiGroup}
\end{dmath}
where $\trace_{\Gamma}\qty(\cdot)$ is partial trace operation over $\Gamma$. We restrict our discussion to the case where the semi-group $\mathcal{T}_t$ is uniformly continuous so that its generator is bounded \citep[Prop 3.8]{fagnola1999quantum}. Let us write $\rho_t \equiv \mathcal{S}_t\qty(\rho)$. The generator of this semi-group is the master equation \citep{fagnola1999quantum} 
\begin{dmath}
\mathcal{L}_{\ast}\qty(\rho_t) \equiv
-i\left[\mathbb{H},\rho_t\right] + \mathbb{L}^{\top}\rho_t\mathbb{L}^{\ast} - \frac{1}{2}\mathbb{L}^{\dagger}\mathbb{L}\rho_t - \frac{1}{2}\rho_t\mathbb{L}^{\dagger}\mathbb{L}. \label{eq:MasterEquation}
\end{dmath}
For our development, we will also require the following definitions:
\begin{defn}\citep{Fagnola2001}
	A density operator $\rho$ is called invariant for a {QDS} $\mathcal{T}_t$ if for all $A \in \mathfrak{B}\qty(\mathscr{H})$, $t\geq 0$, $\trace\qty(\rho \mathcal{T}_t\qty(A))  = \trace\qty(\rho A)$.
\end{defn}
We recall that the set of trace-class Hilbert space operators $\mathfrak{I}_1\qty(\mathscr{H})$ with norm $\norm{\cdot}_1$ is a Banach space \citep[Prop 9.13]{Parthasarathy1992}.
Using the metric induced by the norm $d\qty(\rho_A,\rho_B) \equiv \norm{\rho_A-\rho_B}_1$, we define a closed ball with center $\rho_{\ast}$ and radius $\epsilon$ by the set
\begin{dmath}
\mathscr{B}_{\epsilon}\qty(\rho_{\ast}) = \qty{\rho \in \mathfrak{S}\qty(\mathscr{H}) : \norm{\rho-\rho_{\ast}}_1\leq \epsilon}. \label{eq:BallCenterAtRhoStarRadiusEpsilon}
\end{dmath}
The normalized version of the distance $d\qty(\rho_A,\rho_B)$ is also known as the Kolmogorov distance in the quantum information community \citep{Fuchs1999,MichaelA.Nielsen2001,IngemarBengtsson2017}; see Figure \ref{fig:QubitBall} for an illustration of $\mathscr{B}_{\epsilon}\qty(\rho_{\ast})$ for a qubit system with center $\rho_\ast = \frac{1}{2}\mqty[1 && 0\\0 && 1]$.

We will also define a neighborhood $\mathscr{N}$ of $\rho_\ast$ as a union of balls \eqref{eq:BallCenterAtRhoStarRadiusEpsilon} with various center points and $\mathscr{B}_{\epsilon}\qty(\rho_{\ast}) \subseteq \mathscr{N}$ for some $\epsilon$.
The following proposition shows a basic fact regarding the closedness of the class of density operators.
\begin{prop}\label{prp:SHisCompleteSubspaceOfI1}
	The class of density operators on the Hilbert space $\mathscr{H}$, $\mathfrak{S}\qty(\mathscr{H})$ is a closed subset of the Banach space $\qty(\mathfrak{I}_1\qty(\mathscr{H}),\norm{\cdot}_1)$.
\end{prop}
\begin{pf}
First we recall that a subset of a complete metric space is closed if and only if it is complete \citep[Prop 6.3.13]{Morris2016}. Therefore, we need to show that every Cauchy sequence of density operators $\qty{\rho_n}$ converges to a density operator $\rho_\ast$. Since $\mathfrak{S}\qty(\mathscr{H}) \subseteq \mathfrak{I}_1\qty(\mathscr{H})$, which is a Banach space with respect to the norm $\norm{\cdot}_1$, then $\qty{\rho_n}$ converges to an element in $\mathfrak{I}_1\qty(\mathscr{H})$; i.e., $\rho_\ast \in \mathfrak{I}_1\qty(\mathscr{H})$. Therefore according to the definition of the class of density operators, it remains to show that $\rho_\ast $ is positive and has unity trace. The limit $\rho_\ast$ is positive, since if it is non-positive, then there exists $\epsilon >0$ such that for all $n$, $0 < \epsilon < \norm{ \rho_n - \rho_\ast}_\infty$. However, as $n \rightarrow \infty$,
\begin{align*}
0 < \epsilon < \norm{ \rho_n - \rho_\ast}_\infty \leq  \norm{ \rho_n - \rho_\ast}_1 = 0,
\end{align*}
which is a contradiction. 
The limit $\rho_\ast$ also has unit trace by the following argument. Since $\rho_n$ converges to $\rho_\ast$, then for any $\epsilon >0$ there is an $n$ such that 
\begin{align*}
1 = \norm{ \rho_n }_1 \leq \norm{ \rho_n - \rho_\ast}_1 + \norm{ \rho_\ast}_1 \leq \epsilon + \norm{\rho_\ast}_1.
\end{align*}
However, we notice also that for any $\epsilon >0$, there is an $N_\epsilon \in \mathbb{N}$, such that for every $n,m \geq N_\epsilon$, $\norm{\rho_n - \rho_m}_1 < \epsilon$. 
Fix $n$. Then we have $\norm{\rho_n}_1 \leq \norm{\rho_n - \rho_{N_\epsilon}}_1 +\norm{\rho_{N_\epsilon}}_1$. Taking the limit as $n$ approaches infinity, we obtain 
\begin{dmath*}
\norm{\rho_{\ast}}_1 \leq  \lim\limits_{n\rightarrow \infty} \norm{\rho_n}_1 + \norm{\rho_n - \rho_{\ast}}_1 
\leq  \lim\limits_{n\rightarrow \infty} \norm{\rho_n - \rho_{N_\epsilon}}_1 +\norm{\rho_{N_\epsilon}}_1 + \norm{\rho_n - \rho_{\ast}}_1 
\leq \epsilon + 1.
\end{dmath*}
Since $\epsilon$ can be chosen arbitrarily, then $\norm{\rho_{\ast}}_1 = 1$. Therefore, $\rho_{\ast}$ is indeed a density operator.\qed
\end{pf}

\section{Stability of Quantum Systems I: A Lyapunov Criterion}
In this section, we will introduce a Lyapunov stability concept for quantum system invariant density operators. Before we define the stability condition, we first show that the set of invariant density operators of the {QDS} \eqref{eq:TSemiGroup} is both closed and convex.
\begin{prop}\label{prp:InvariantDensityOperatorConvexClosed}
	\textcolor{black}{If the set of invariant density operators $\mathscr{C}_\ast$ for a QDS of the form \eqref{eq:TSemiGroup} is non-empty}, then it is convex and closed in $\qty(\mathfrak{I}_1\qty(\mathscr{H}),\norm{\cdot}_1)$. 
\end{prop}
\begin{pf}
	The convexity of $\mathscr{C}_\ast$ follows directly from the fact that for any $\rho_1,\rho_2 \in \mathscr{C}_\ast$, then for any $\lambda \in \qty[0,1]$, $\trace\qty(\qty(\mathcal{S}_t\qty(\lambda \rho_1 + (1-\lambda)\rho_2) - \qty(\lambda \rho_1 + (1-\lambda)\rho_2))A) =  0$, for all $A \in \mathfrak{B}\qty(\mathscr{H})$, and $t \geq 0$.
	Notice that since $\mathscr{C}_\ast$  is convex, the closedness of  $\mathscr{C}_\ast$ on $\qty(\mathfrak{I}_1\qty(\mathscr{H}),\norm{\cdot}_1)$ is equivalent to the closedness of $\mathscr{C}_\ast$  in the weak topology  \citep[Thm III.1.4 ]{Conway1985}.
	To show that $\mathscr{C}_\ast$ is closed in the weak topology, suppose that $\qty{\rho_n}$ is a net in $\mathscr{C}_\ast$ that converges weakly to $\rho_\ast$, $\rho_n \xrightharpoondown{w} \rho_\ast$. Then we have to show that $\rho_\ast \in \mathscr{C}_\ast$.  We observe that the linearity of the semi-group $\mathcal{S}_t$ implies that for all $A \in \mathfrak{B}\qty(\mathscr{H})$, $\rho_n$ in the net $\qty{\rho_n}$, and $t\geq 0$
	\begin{dmath*}
	\trace\qty(\qty(\mathcal{S}_t\qty(\rho_\ast)-\rho_\ast)A) = \trace\qty(\qty(\mathcal{S}_t\qty(\rho_\ast - \rho_n))A) + \trace\qty(\qty(\mathcal{S}_t\qty(\rho_n) - \rho_n)A) + \trace\qty(\qty(\rho_n - \rho_\ast)A)\\ = \trace\qty(\qty(\mathcal{S}_t\qty(\rho_\ast - \rho_n))A)  + \trace\qty(\qty(\rho_n - \rho_\ast)A) = 
	\trace\qty(\qty(\rho_\ast - \rho_n)\qty(\mathcal{T}_t\qty(A) - A)).
	\end{dmath*}
Since $\mathcal{T}_t\qty(A) \in \mathfrak{B}\qty(\mathscr{H})$, then for any $\epsilon >0$, there exists a $\rho_m \in \qty{\rho_n}$ such that $\abs{\trace((\rho_\ast - \rho_m)(\mathcal{T}_t\qty(A) - A))} < \epsilon$.
However,  $\epsilon>0$ can be selected arbitrarily. Therefore, $\rho_\ast  \in \mathscr{C}_\ast$.\qed
\end{pf}
The convexity of $\mathscr{C}_\ast$ has been previously mentioned in \citep{Schirmer2010} for the finite dimensional case.  The last proposition also tells us that it is irrelevant to ask about the convergence of a density operator evolution to a particular point inside $\mathscr{C}_\ast$ when it consists more than one density operator, since every neighborhood $\mathscr{N}$ of $\rho_\ast\in\mathscr{C}_\ast$ must contain other invariant density operators for which a trajectory can eventually converge to. 
In what follows, we will examine the convergence to set of invariant density operators in a Banach space $\qty(\mathfrak{I}_1\qty(\mathscr{H}),\norm{\cdot}_1)$. The distance between a point $\sigma \in \mathfrak{S}\qty(\mathscr{H})$ and the closed convex set $\mathscr{C}_\ast \subseteq \mathfrak{S}\qty(\mathscr{H})$ can  be naturally defined by
\begin{dmath}
d\qty(\sigma,\mathscr{C}_\ast) = \inf_{\rho \in \mathscr{C}_\ast} \norm{\sigma-\rho}_1. \label{eq:dtoC_ast}
\end{dmath}
We define the following stability notions:
\begin{defn}\label{def:StabilityOfDensityOperator}
	Let $\mathscr{C}_{\ast} \subset \mathfrak{S}\qty(\mathscr{H}) $ be a \textcolor{black}{non-empty} convex set of invariant density operators of a QDS \eqref{eq:TSemiGroup}. Suppose that $\mathscr{N} \subset \mathfrak{S}\qty(\mathscr{H})$, where $\mathscr{C}_{\ast}$ is a strict subset of $\mathscr{N}$, and the system is initially at density operator $\rho \in \mathscr{N}$. Then, we say the closed convex set of invariant density operators $\mathscr{C}_{\ast}$ is,
	\begin{enumerate}
		\item Lyapunov stable if for every $\varepsilon >0$,  there exists $\delta(\varepsilon) >0$ such that $d\qty(\rho,\mathscr{C}_{\ast}) < \delta(\varepsilon)$ implies $d\qty(\rho_t,\mathscr{C}_{\ast}) < \varepsilon$ for all $t\geq 0$.
		\item Locally asymptotically stable, if it is Lyapunov stable, and there exists $\delta >0$, such that $d\qty(\rho,\mathscr{C}_{\ast}) < \delta$ implies $\lim\limits_{t\rightarrow\infty} d\qty(\rho_t,\mathscr{C}_{\ast}) = 0$.
		\item Locally exponentially stable, if there exists $\beta,\gamma,\delta >0$  such that $d\qty(\rho,\mathscr{C}_{\ast}) < \delta$ implies, $d\qty(\rho_t,\mathscr{C}_{\ast}) \leq \beta d\qty(\rho,\mathscr{C}_{\ast}) \exp\qty(-\gamma t)$ for all $t\geq 0$.
	\end{enumerate}
	If $\mathscr{N} = \mathfrak{S}\qty(\mathscr{H})$, such that $\delta$ can be chosen arbitrarily in 2) and 3), we say $\mathscr{C}_{\ast}$ is a globally asymptotically, or exponentially stable respectively.
\end{defn}
Before we prove our first main result, we need to establish the following facts. Let $A \in \mathfrak{B}\qty(\mathscr{H})$ be a self-adjoint operator where the spectrum of $A$ is strictly increasing, that is $A = \sum_{i=0}^{\infty}p_i P_i$, where $\qty{P_i}$ is a set of orthogonal projection operators (not necessarily of rank one), and the spectrum values of $A$ satisfy $p_0 < p_1 < p_2 < \cdots < p_n < \cdots$. The set of \emph{ground states} of $A$, is given by $\mathscr{W}_A = \qty{\rho \in \mathfrak{S}\qty(\mathscr{H}): \trace(\rho A) = p_0}$ \citep{Pan2016}.
\begin{lem}\label{lem:SuperPosition}
Let two unit vectors $\ket{\psi_i} \in \mathscr{H}, i=1,2$, $\bra{\psi_2}\ket{\psi_1}=0$ be eigen-vectors of a self-adjoint operator $A$ corresponding to the same spectrum value $p$; i.e., $A\ket{\psi_i} = p\ket{\psi_i}, i=1,2$. Then any superposition vector of $\ket{\psi_1}$ and $\ket{\psi_2}$ is also an eigen-vector of $A$ corresponding to $p$. In particular, if $\ket{\psi_i}\bra{\psi_i}, i=1,2$ are ground states of $A$, then if $\ket{\psi}$ is a superposition of $\ket{\psi_1}$ and $\ket{\psi_2}$, then $\ket{\psi}\bra{\psi}$ is a also ground state of $A$.
\end{lem}
\begin{pf}
	Immediate; see also \citep[p. 71]{GantMacher1959}.\qed
\end{pf}
\begin{lem}\label{lem:ConstantAcrossCStar}
	Suppose there exists a self-adjoint operator $A \in \mathfrak{B}\qty(\mathscr{H})$ where the spectrum of $A$ is strictly increasing such that for a \textcolor{black}{non-empty} closed convex set of density operators $\mathscr{C}_\ast$, and a neighborhood $\mathscr{N}$ of $\mathscr{C}_\ast$, 
\begin{align}
\inf_{\rho_\ast \in \mathscr{C}_\ast}\trace\qty(A\qty(\rho - \rho_{\ast})) > 0, \; \forall \rho \in \mathscr{N}\backslash \mathscr{C}_\ast. \label{eq:StrictMinimaOnNeighbourhood}
\end{align}
Then $\trace\qty(A\qty(\rho_\ast))$ is constant for all  $\rho_\ast \in \mathscr{C}_\ast$. Moreover, $\mathscr{C}_\ast=\mathscr{W}_A$.
\end{lem}
\begin{pf}
	Let the operator $A = \sum_{i=0}^{\infty}p_i P_i$, where the spectrum values of $A$ satisfy $p_0 < p_1 < p_2 < \cdots < p_n < \cdots$.
	Suppose $\sup_{\rho_\ast \in \mathscr{C}_\ast} \trace\qty(A\rho_\ast) = \bar{\alpha}$, and $\inf_{\rho_\ast \in \mathscr{C}_\ast} \trace\qty(A\rho_\ast) = \underset{\bar{}}{\alpha}$,  where $\underset{\bar{}}{\alpha} < \bar{\alpha}$. 
	\textcolor{black}{Consider the topology on $\mathfrak{S}(\mathscr{H})$ induced by norm $\norm{\cdot}_1$. By definition, $\mathscr{C}_\ast \subset \mathfrak{S}(\mathscr{H})$, and from Propositions \ref{prp:SHisCompleteSubspaceOfI1} and \ref{prp:InvariantDensityOperatorConvexClosed}, $\mathscr{C}_\ast$ is not a clopen set (a set that is both closed and open at the same time) in $\mathfrak{S}(\mathscr{H})$. Hence,  it is clear that the boundary of $\mathscr{C}_\ast$ is non-empty in $\mathfrak{S}(\mathscr{H})$; see \cite[103, exercise 7]{Mendelson2003}. Therefore, } 
	pick any $\rho_{\ast\ast}$ at the boundary surface of $\mathscr{C}_\ast$ and choose $\rho' \in \mathscr{N} \backslash \mathscr{C}_\ast$ such that for any $\lambda \in (0,1]$, $\lambda \rho' + (1-\lambda) \rho_{\ast\ast}$ does not belong to $\mathscr{C}_\ast$. Let $\delta = \bar{\alpha} - \underset{\bar{}}{\alpha}$, $\trace(A\rho') = \bar{\alpha} + \varepsilon_0$, $\trace(A\rho_{\ast\ast}) = \underset{\bar{}}{\alpha} + \varepsilon_1$ where $\varepsilon_0>0$, and $\varepsilon_1 \in [0,\delta]$.
	Therefore, setting $\lambda \leq (\delta - \varepsilon_1)/(\delta+\varepsilon_0 - \varepsilon_1)$, we obtain
	\begin{align}
	\trace(A(\lambda \rho' + (1-\lambda) \rho_{\ast\ast})) &\leq \bar{\alpha},
	\end{align}
	which violates the inequality \eqref{eq:StrictMinimaOnNeighbourhood}.\\
	Now we will show that $\mathscr{C}_\ast \subseteq \mathscr{W}_A$, which means that $\trace\qty(A\rho_\ast) = p_0$, for any $\rho_\ast \in  \mathscr{C}_\ast$.  To see this, suppose  $\mathscr{C}_\ast \nsubseteq \mathscr{W}_A$. Then, there exists a $\rho_\ast \in  \mathscr{C}_\ast$ such that $\trace\qty(A\rho_\ast) = \hat{p} > p_0$. But since $\trace\qty(A\rho_\ast)$ has to be constant across all  $\rho_\ast \in \mathscr{C}_\ast$, then $\trace\qty(A\rho_\ast) = \hat{p} > p_0$ for all $\rho_\ast \in \mathscr{C}_\ast$. Select $\rho' \in \mathscr{W}_A \backslash \mathscr{C}_\ast$. Then for any $\rho_\ast \in \mathscr{C}_\ast$, there exists $\lambda' \in (0,1)$ such that , $\rho = \lambda' \rho_\ast + (1-\lambda')\rho' \in \mathscr{N}\backslash \mathscr{C}_\ast$. However, if this is the case, then $\trace(A\rho) = (1-\lambda') p_0 + \lambda' \hat{p} \leq \hat{p}$ which contradicts \eqref{eq:StrictMinimaOnNeighbourhood}. 
	\\
	Next since $\mathscr{C}_\ast \subseteq \mathscr{W}_A$, suppose there exists $\rho' \in  \mathscr{W}_A \backslash \mathscr{C}_\ast$. Then for an arbitrary $\rho_\ast \in \mathscr{C}_\ast$, there also exists a $\lambda' \in (0,1)$ such that, $\rho = \lambda' \rho_\ast + (1-\lambda')\rho' \in \mathscr{N}\backslash \mathscr{C}_\ast$. However, if this is the case, then $\trace(A\rho) = p_0$, which contradicts \eqref{eq:StrictMinimaOnNeighbourhood}.\qed
\end{pf}

\begin{lem}\label{lem:LowerBound}
	Suppose the condition in Lemma \ref{lem:ConstantAcrossCStar} is satisfied. Then there exists $\kappa>0$ such that
	\begin{align*}
	\kappa d\qty(\rho,\mathscr{C}_\ast)^2 \leq& \inf_{\rho_\ast \in \mathscr{C}_\ast}\trace\qty(A\qty(\rho - \rho_{\ast})),
	\end{align*}
	for all $\rho \in \mathscr{N} \backslash \mathscr{C}_\ast.$
\end{lem}
\begin{pf}
	Let the operator $A$ be as in the proof of Lemma \ref{lem:ConstantAcrossCStar}. By Lemma \ref{lem:ConstantAcrossCStar}, $\mathscr{C}_\ast=\mathscr{W}_A$. Without loss of generality, we can assume that $p_0 = 0$ (otherwise we can consider $A' = A - p_o \mathbb{1}$) and there exists a collection of $r_p$ unit vectors $\qty{\ket{j}}$ such that $\bra{i}\ket{j}=\delta_{ij}$, and $\bra{j}A\ket{j}=p_0$ for all $j = 0,\cdots, r_p-1$. Without loss of generality, we can assume that $\ket{j}$ is a number vector, where $\ket{0} = (1,0,\cdots);\ket{1} = (0,1,\cdots)$; etc. Otherwise there is a unitary transformation $W$ such that $A' = W^\dagger A W$ is diagonal. Therefore, by Lemma \ref{lem:SuperPosition}, any invariant density operator $\rho_\ast$ is a convex combination of unit vectors $\qty{\ket{s_j}}$ that are superpositions of $\qty{\ket{j}}$.
	\\
	Let us consider a unitary vector $\ket{\psi_i} = \qty(a_{i0},a_{i1},\cdots)$, where $a_{ij} \in \mathbb{C}$, and $\sum_{j=0}^{\infty}\abs{a_{ij}}^2 = 1$. Select $\rho = \sum_{i} \lambda_i \ket{\psi_i}\bra{\psi_i}$, where $\lambda_i >0$, $\bra{\psi_i}\ket{\psi_j} = \delta_{ij}$ and $\sum_{i} \lambda_i = 1$. Notice that, we can also write any $\rho_\ast \in \mathscr{C}_\ast$ in the form $\rho_\ast = \sum_{i} \lambda_i \hat{\rho}_i$, where $\qty{\hat{\rho}_i} \in \mathscr{C}_\ast$.
	\\
	Using the convexity of the norm $\norm{ \cdot }_1$, we observe that
	\begin{dmath*}
	d(\rho,\mathscr{C}_\ast) = \inf_{\rho_\ast \in \mathscr{C}_\ast} d(\rho,\rho_\ast) = \inf_{\qty{\hat{\rho}_i} \in \mathscr{C}_\ast} d\qty(\sum_{i} \lambda_i \ket{\psi_i}\bra{\psi_i},\sum_{i} \lambda_i \hat{\rho}_i) 
	=
	\inf_{\qty{\hat{\rho}_i} \in \mathscr{C}_\ast} \norm{\sum_{i} \lambda_i \left(\ket{\psi_i}\bra{\psi_i} - \hat{\rho}_i\right) }_1\\
	=
	\inf_{\qty{\hat{\rho}_i} \in \mathscr{W}_A} \norm{\sum_{i} \lambda_i \left(\ket{\psi_i}\bra{\psi_i} - \hat{\rho}_i\right) }_1
	\leq \inf_{\qty{\hat{\rho}_i} \in \mathscr{W}_A} \sum_{i} \lambda_i\norm{ \ket{\psi_i}\bra{\psi_i} - \hat{\rho}_i}_1\\
	= \sum_{i} \lambda_i \inf_{\hat{\rho}_i \in \mathscr{W}_A} \norm{ \ket{\psi_i}\bra{\psi_i} - \hat{\rho}_i}_1
	=  \sum_{i} \lambda_i d\qty(\ket{\psi_i}\bra{\psi_i},\mathscr{W}_A).
	\end{dmath*}
	Moreover, due to the convexity of the quadratic function, we further obtain
	\begin{dmath*}
	d(\rho,\mathscr{C}_\ast)^2 \leq \left(\sum_{i} \lambda_i d\qty(\ket{\psi_i}\bra{\psi_i},\mathscr{W}_A)\right)^2
	\leq \sum_{i} \lambda_i d\qty(\ket{\psi_i}\bra{\psi_i},\mathscr{W}_A)^2.
	\end{dmath*}
	Now, for any $\ket{\psi}$ and $\ket{\phi} \in \mathscr{H}$, we have \citep[\textsection 9.2.3]{MichaelA.Nielsen2001}
	\begin{align*}
	\norm{ \ket{\psi}\bra{\psi} - \ket{\phi}\bra{\phi}}_1 =& 2 \sqrt{1 - \abs{\bra{\psi}\ket{\phi}}^2} .
	\end{align*}
	Therefore
	\begin{dmath*}
	d(\ket{\psi_i}\bra{\psi_i},\mathscr{W}_A )^2 = \qty(\inf_{\hat{\rho}_i \in \mathscr{W}_A} \norm{ \ket{\psi_i}\bra{\psi_i} - \hat{\rho}_i}_1)^2 \leq 4 \qty(1 - \sum_{j=0}^{r_p-1} \abs{\bra{\psi_i}\ket{j}}^2) = 4 \sum_{j=r_p} \abs{a_{ij}}^2.
	\end{dmath*}	
	Now we also observe that 
	\begin{align*}
	\trace\qty(A (\rho-\rho_\ast)) = \sum_{i} \lambda_i \sum_{j=r_p}p_j \abs{a_{ij}}^2, \; j > 0.
	\end{align*}
	Hence
	\begin{align*}
	\dfrac{\trace\qty(A (\rho-\rho_\ast))}{d\qty(\rho,\mathscr{C}_\ast)^2} \geq \dfrac{\sum_{i} \lambda_i \sum_{j=r_p}p_j \abs{a_{ij}}^2}{4 \sum_{i} \lambda_i \sum_{j=r_p} \abs{a_{ij}}^2} 
	\geq {\dfrac{p_1}{4}>0}.
	\end{align*}
	This completes the proof.\qed
\end{pf}

\begin{cor} \label{cor:WeakConverge}
Let $\mathscr{C}_\ast$ be a closed convex set of invariant density operators and $\mathscr{N}$ be a neighborhood of $\mathscr{C}_\ast$. Suppose $\qty{\rho_n}$ is a sequence of density operators in $\mathscr{N} \backslash \mathscr{C}_\ast$. If there exists a self-adjoint operator $A \in \mathfrak{B}\qty(\mathscr{H})$ satisfying the condition in Lemma \ref{lem:ConstantAcrossCStar} and $\lim\limits_{n\rightarrow \infty}\trace\qty(A\qty(\rho_n - \rho_\ast) ) = 0$ for a $\rho_\ast \in \mathscr{C}_\ast$, then  $\lim\limits_{n\rightarrow \infty} d\qty(\rho_n,\mathscr{C}_\ast) = 0$.
\end{cor}
\begin{pf}
Let$\lim\limits_{n\rightarrow \infty}\trace\qty(A\qty(\rho_n - \rho_\ast) ) = 0$ for some $A \in \mathfrak{B}\qty(\mathscr{H})$ satisfying Lemma \ref{lem:ConstantAcrossCStar}. 
By Lemma \ref{lem:LowerBound} there exists a $\kappa > 0$, such that $\kappa d\qty(\rho,\mathscr{C}_\ast)^2 \leq \inf_{\rho_\ast \in \mathscr{C}_\ast}\trace\qty(A\qty(\rho_n - \rho_{\ast}))$. Therefore as $n$ goes to infinity, $\lim\limits_{n\rightarrow\infty} d\qty(\rho_n,\mathscr{C}_\ast) = 0$, which completes the proof.\qed
\end{pf}
The following theorem is the first main result of this article, which relates the stability notions in Definition \ref{def:StabilityOfDensityOperator} to inequalities for the generator of a candidate Lyapunov operator. 
\\
As a reminder, we notice that by \eqref{eq:LandGrelation} we have the following equality $\Paverage{\mathcal{G}\qty(V_t)} =
\trace\qty(\mathcal{L}\qty(V)\rho_t) $. 
\begin{thm}\label{thm:LyapunovStability}
	Let $V \in \mathfrak{B}\qty(\mathscr{H})$ be a self-adjoint operator with strictly increasing spectrum values such that
	\begin{dmath}
	{\inf_{\rho_\ast \in \mathscr{C}_\ast}\trace\qty(V\qty(\rho - \rho_{\ast})) > 0}, {\forall \rho \in \mathfrak{S}\qty(\mathscr{H}) \backslash  \mathscr{C}_\ast}.
	\label{eq:V_rho_ast_0}
	\end{dmath}

	Using the notation of Definition \ref{def:StabilityOfDensityOperator}, the following properties hold:
	\begin{enumerate}
		\item If
		\begin{dmath}
		\trace\qty(\mathcal{L}\qty(V)\rho) \leq 0, \; {\forall \rho \in \mathfrak{S}\qty(\mathscr{H})\backslash \mathscr{C}_\ast}, \label{eq:LyapunovStability}
		\end{dmath}
		then $\mathscr{C}_\ast$ is Lyapunov stable.
		
		\item If 
		\begin{dmath}
		\trace\qty(\mathcal{L}\qty(V)\rho) < 0, \; {\forall \rho \in \mathfrak{S}\qty(\mathscr{H})\backslash \mathscr{C}_\ast},
		\label{eq:AsymptoticallyStableCondition}
		\end{dmath}
		then $\mathscr{C}_\ast$ is globally asymptotically stable.
		
		\item If there exists $\gamma >0$ and $\zeta \in \mathbb{R}$ such that
		\begin{dmath}
		\trace\qty(\mathcal{L}\qty(V)\rho) \leq -\gamma \trace\qty(V\rho) + \zeta <0  \; {\forall \rho \in \mathfrak{S}\qty(\mathscr{H}) \backslash \mathscr{C}_\ast}, \label{eq:ExponentiallyStableCondition}
		\end{dmath}
		then $\mathscr{C}_\ast$ is globally exponentially stable.
	\end{enumerate}
\end{thm}
\begin{pf}
Let $\mathscr{N} = \mathfrak{S}\qty(\mathscr{H})$. Let us begin by proving the first part. Suppose $\varepsilon >0$ is selected. Then, we can take $\varepsilon' \in (0,\varepsilon]$ such that $\mathscr{B}_{\varepsilon'}\qty(\mathscr{C}_{\ast}) \subseteq \mathscr{N}$. Observe that by \eqref{eq:V_rho_ast_0} and Lemma \ref{lem:LowerBound}, there exists a $\kappa >0$ such that for any $\rho \in  \mathscr{B}_{\varepsilon'}\qty(\mathscr{C}_{\ast}) $
\begin{dmath}
\kappa d\qty(\rho,\mathscr{C}_\ast)^2 \leq \inf_{\rho_\ast \in \mathscr{C}_\ast} \trace\qty(V\qty(\rho-\rho_\ast)). \label{eq:NBetaWeaklyCompact}
\end{dmath}
Let $V_\ast = \sup_{\rho_\ast \in \mathscr{C}_\ast} \trace\qty(V\rho_\ast)$.
Therefore if we select $V_\ast < \beta < V_\ast + \kappa (\epsilon')^2$, then the set $\mathscr{N}_{\beta} = \qty{ \rho \in \mathscr{N} : \trace\qty(V\rho)\leq \beta}$
is a strict subset of $\mathscr{B}_{\varepsilon'}\qty(\rho_{\ast})$.
Furthermore, since $\beta > V_\ast$ and $\trace\qty(V\rho) \leq V_\ast + \norm{V}_\infty d\qty(\rho,\mathscr{C}_\ast)$, selecting $\delta < \qty(\beta - V_\ast)/\norm{V}_\infty$ implies $\trace\qty(V\rho) < \beta$ for all $\rho \in \mathscr{B}_{\delta}\qty(\mathscr{C}_{\ast})$.
Therefore, we have the following relation
\begin{dmath*}
{\mathscr{B}_{\delta}\qty(\mathscr{C}_{\ast}) \subset \mathscr{N}_{\beta} \subset \mathscr{B}_{\varepsilon'}\qty(\mathscr{C}_{\ast})}.
\end{dmath*}
Hence, $\rho \in \mathscr{B}_{\delta}\qty(\mathscr{C}_{\ast})$ implies $\rho \in \mathscr{N}_{\beta}$.  Since $\trace\qty(\mathcal{L}\qty(V)\rho) \leq 0$ for all $\rho \in \mathscr{N}\backslash \mathscr{C}_\ast$, if  system density operator $\rho$ is initially in $\mathscr{N}_{\beta}$, then the expected value of operator $V$ will be non-increasing, $
{\trace\qty(V\rho_t) \leq \trace\qty(V\rho) \leq \beta, \forall t \geq 0}.$
This implies that $\rho_t \in \mathscr{N}_{\beta}, \forall t \geq 0$, which shows $\rho_t \in \mathscr{B}_{\varepsilon'}\qty(\mathscr{C}_{\ast})$. Furthermore, this last implication implies that if initially $d\qty(\rho,\mathscr{C}_{\ast}) < \delta(\varepsilon)$, then $d\qty(\rho_t,\mathscr{C}_{\ast}) < \varepsilon$ for all $t \geq 0$.
\\
For the second part, using the same argument as in the previous part, we may choose $\delta>0$ such that initially $\rho \in \mathscr{B}_{\delta}\qty(\mathscr{C}_\ast) \subset \mathscr{N}_\beta \subset \mathscr{N}$, for some $\beta > V_\ast$. Therefore, since  $\trace\qty(\mathcal{L}\qty(V)\rho) < 0$ for all $\rho \in \mathscr{N}\backslash \mathscr{C}_\ast$, $\trace\qty(V\rho_t)$ is monotonically decreasing. Therefore $\trace\qty(V\rho_t)<\trace\qty(V\rho_s) < \trace\qty(V\rho) < \beta$ for any $0 < s < t$. Hence $\rho_s,\rho_t$ also belong to $\mathscr{N}_\beta$. This implies that there exists a sequence of density operators $\qty{\rho_n}$, such that $\trace\qty(V\qty(\rho_m - \rho_n)) <0$ for any $m>n$, where $\rho_n \equiv \mathcal{S}_{t_n}\qty(\rho)$ and $0 \leq t_0 <t_1 < \cdots < t_n$, $t_n \rightarrow \infty$ as $n \rightarrow \infty$. Since the spectrum of $V$ is non-decreasing, $\trace\qty(V\rho_n)$ is lower bounded. Hence, there exists a $\rho_c \in \mathscr{N}_\beta$ such that $\lim\limits_{n\rightarrow \infty} \trace\qty(V\qty(\rho_n - \rho_c)) = 0$ and $\trace\qty(V\qty(\rho_c - \rho_n)) <0$ for all $n$. Suppose $\rho_c \notin \mathscr{C}_\ast$. Then for any $s>0$, $\trace\qty(V\mathcal{S}_s\qty(\rho_c)) < \trace\qty(V\rho_c)$. Therefore, there exists an $n$ such that $\trace\qty(V\qty(\rho_c - \rho_n)) >0$, which is a contradiction. Therefore, $\rho_c \in \mathscr{C}_\ast$. Corollary \ref{cor:WeakConverge} and \eqref{eq:dtoC_ast} then imply that $\lim\limits_{n\rightarrow \infty} d\qty(\rho_n, \mathscr{C}_\ast) \leq \lim\limits_{n\rightarrow \infty} \norm{\rho_n - \rho_c}_1 = 0$.
\\    
For the global exponentially stable condition, the previous part shows that the negativity of $\trace\qty(\mathcal{L}\qty(V)\rho)$ for all $\rho \in \mathscr{N}\backslash \mathscr{C}_\ast$ implies the existence of a $\rho_c \in \mathscr{C}_\ast$ such that $\lim\limits_{t\rightarrow \infty} \norm{\rho_t-\rho_c }_1=0$. Using the First Fundamental Lemma of Quantum Stochastic Calculus \citep[Prop 25.1,Prop 26.6]{Parthasarathy1992} to switch the order of the integration and quantum expectation, we obtain
\begin{dmath*}
\trace\qty(V\rho_t) - \trace\qty(V\rho_s)  = \Paverage{V_t} - \Paverage{V_s} =  \Paverage{\int_{s}^{t} \mathcal{G}\qty(V_{\tau}) d\tau}  = \int_{s}^{t} \Paverage { \mathcal{G}\qty(V_{\tau})} d\tau.
\end{dmath*}
Therefore, by \eqref{eq:ExponentiallyStableCondition} we obtain 
\begin{dmath*}
\trace\qty(V\rho_t) \leq \qty(\trace\qty(V\rho_s) - \dfrac{\zeta}{\gamma})e^{\gamma(s-t)} + \dfrac{\zeta}{\gamma}.
\end{dmath*}
Taking $s=0$, there exists $\kappa >0$ such that $\kappa d(\rho_t,\mathscr{C}_\ast)^2 \leq \kappa \norm{\rho_t - \rho_c}_1^2 \leq \trace\qty(V(\rho_t - \rho_c)) \leq \trace(V\rho_t) - \frac{\zeta}{\gamma} \leq \qty(\trace(V\rho) - \frac{\zeta}{\gamma})e^{- \gamma t}$. Consequently, we obtain
\begin{dmath*}
k d(\rho_t,\mathscr{C}_\ast)^2  \leq  \qty(\trace(V\rho) -\frac{\zeta}{\gamma})e^{-\gamma t},
\end{dmath*}
which completes the proof.\qed 	 	
\end{pf}
Theorem \ref{thm:LyapunovStability} establishes sufficient conditions for the stability notions given in Definition \ref{def:StabilityOfDensityOperator} in global settings. Observe that we cannot make a localized version of Theorem \ref{thm:LyapunovStability} by interchanging $\mathfrak{S}\qty(\mathscr{H})$ in \eqref{eq:V_rho_ast_0} with a smaller neighborhood $\mathscr{N}$ of $\mathscr{C}_\ast$. This is the consequence of Lemma \ref{lem:ConstantAcrossCStar},  that states $\mathscr{C}_\ast = \mathscr{W}_V$ when \eqref{eq:V_rho_ast_0} is satisfied for a neighborhood $\mathscr{N}$ of $\mathscr{C}_\ast$. Therefore \eqref{eq:V_rho_ast_0} will hold globally once it holds for a neighborhood $\mathscr{N}$ of $\mathscr{C}_\ast$, which leads to global stability condition.
\\
We have seen in Proposition \ref{prp:InvariantDensityOperatorConvexClosed} that the convexity of the set of the invariant density operators prevents the possibility of multiple isolated invariant density operators.
\\
There is still a question whether in a uniformly continuous QDS, there \emph{is} a local asymptotic stability phenomenon that does not hold globally? In the following theorem, we will see that the answer is negative.
\begin{figure}
	\centering
	\includegraphics[width=0.5\textwidth]{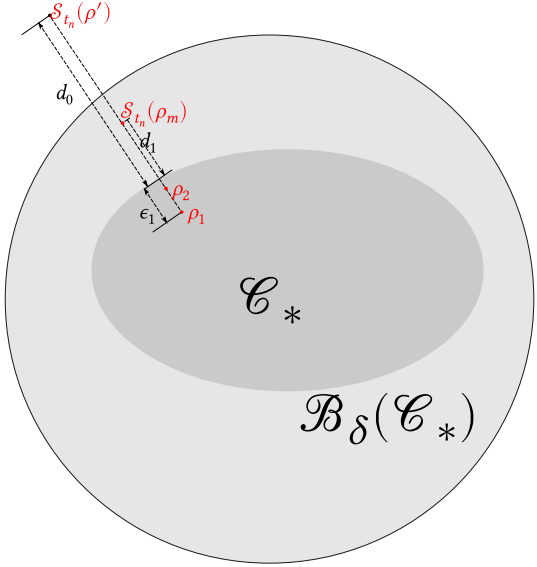}
	\caption{Illustration of the set of invariant density operators %
		in the context of the proof of Theorem \ref{thm:Global_Local_Asymptotic}}.
	\label{fig:Global_Local_Asymptotic}	
\end{figure}

\begin{thm}\label{thm:Global_Local_Asymptotic}
	The set of invariant density operators of a QDS \eqref{eq:TSemiGroup} is globally asymptotically stable if and only if it is locally asymptotically stable.
\end{thm}
\begin{pf}
It is obvious that the global asymptotic stability implies local asymptotic stability. Now suppose that the set of invariant density operators $\mathscr{C}_\ast$ is locally asymptotically stable. Assume that the stability of $\mathscr{C}_\ast$ does not hold globally. Then there exists $\rho' \in \mathfrak{S}\qty(\mathscr{H})$, such that $d\qty(\mathcal{S}_t\qty(\rho'),\mathscr{C}_\ast)$ does not converge to zero as $t \rightarrow \infty$. This means that there exists an $\epsilon' >0$, such that for any $t\geq 0$, there is a $t_{\epsilon'} \geq t$ such that $d\qty(\mathcal{S}_{t_{\epsilon'}}\qty(\rho'),\mathscr{C}_\ast) \geq \epsilon'$. This is equivalent to the existence of a sequence of times $\qty{t_n}, 0 \leq t_0 < t_1 \cdots <t_n$, $t_n \rightarrow \infty$, such that for any $n$, $d\qty(\mathcal{S}_{t_n}\qty(\rho'),\mathscr{C}_\ast) \geq \epsilon'$.
\\
Now, let us examine the following convex set
\begin{align}
\mathscr{W} =& \left\{ \rho\in \mathfrak{S}\qty(\mathscr{H}): \rho = \lambda \rho' + \qty(1-\lambda)\rho_\ast, \rho_\ast \in \mathscr{C}_\ast, \right.\nonumber\\
&\left.\lambda \in \qty(0,1) \right\}.
\end{align}
For any $\delta > 0$, $\mathscr{B}_{\delta}\qty(\mathscr{C}_\ast)\cap \mathscr{W} \neq \varnothing$. The local asymptotic stability property of $\mathscr{C}_\ast$ implies that for some $\delta >0$, any $\rho_m \in \mathscr{B}_{\delta}\qty(\mathscr{C}_\ast)\cap \mathscr{W}$ will converge to $\mathscr{C}_\ast$; i.e., for all $\epsilon >0$ there exists a $t_{\epsilon}$ such that for all $t \geq t_{\epsilon}$, $d\qty(\mathcal{S}_t\qty(\rho_m),\mathscr{C}_{\ast}) < \epsilon$.  Also, there exists a $0<\bar{\lambda}<1$ such that for any $\rho_\ast \in \mathscr{C}_\ast$, $\qty(\bar{\lambda} \rho' + \qty(1-\bar{\lambda}) \rho_\ast)$ belongs to $\mathscr{B}_{\delta}\qty(\mathscr{C}_\ast)\cap \mathscr{W}$. 
Now we can choose $\epsilon = \qty(\bar{\lambda} \epsilon')/2$, and select $t_n \geq t_\epsilon$. Therefore, for this $t_n$ we have $d\qty(\mathcal{S}_{t_n}\qty(\rho'),\mathscr{C}_\ast) \geq \epsilon'$ and
\begin{align}
d\qty(\mathcal{S}_{t_n}\qty(\bar{\lambda} \rho'+\qty(1-\bar{\lambda})\rho_\ast),\mathscr{C}_{\ast}) < \qty(\bar{\lambda} \epsilon')/2, & \forall \rho_\ast \in \mathscr{C}_\ast. \label{eq:d_rho_m_Less_Than}
\end{align} 
For the sake of simplicity let us define for any $\rho_e$ belonging to $\mathscr{C}_\ast$:
\begin{align}
d_0 \equiv& d\qty(\mathcal{S}_{t_n}\qty(\rho'),\mathscr{C}_\ast), \label{eq:d_0}\\
d_1\qty(\rho_e) \equiv& d\qty(\bar{\lambda}  \mathcal{S}_{t_n}\qty(\rho') + \qty(1-\bar{\lambda}) \rho_e,\mathscr{C}_{\ast}) < \qty(\bar{\lambda} \epsilon')/2, \label{eq:d_1}%
\end{align}
where we use the inequality \eqref{eq:d_rho_m_Less_Than} to obtain \eqref{eq:d_1}. 
Next, we claim that there exists a $\rho_e \in \mathscr{C}_\ast$ such that:
\begin{align}
\bar{\lambda} d_0 \leq& d_1\qty(\rho_e) + \frac{\bar{\lambda}\epsilon'}{2}.\label{eq:d_Clambda_1}
\end{align}
To see this, let $\epsilon_1 = \frac{\bar{\lambda}\epsilon'}{4(1-\bar{\lambda})}$ and $\epsilon_2 = \frac{\bar{\lambda}\epsilon'}{4}$. By definition, there exist density operators $\rho_1,\rho_2 \in \mathscr{C}_\ast$ such that %
\begin{align}
d\qty(\mathcal{S}_{t_n}\qty(\rho'),\rho_1) \leq& d_0 + \epsilon_1,\\
d\qty(\bar{\lambda}\mathcal{S}_{t_n}\qty(\rho') + \qty(1-\bar{\lambda}) \rho_1,\rho_2) \leq& d_1\qty(\rho_1) + \epsilon_2.
\end{align}
Therefore, we can write
\begin{align*}
d_0 \leq& \norm{\mathcal{S}_{t_n}\qty(\rho') - \rho_2}_1\\
\leq& \norm{\mathcal{S}_{t_n}\qty(\rho') - \qty(\bar{\lambda}\mathcal{S}_{t_n}\qty(\rho') + \qty(1-\bar{\lambda}) \rho_1)}_1 \\
+& \norm{\qty(\bar{\lambda}\mathcal{S}_{t_n}\qty(\rho') + \qty(1-\bar{\lambda}) \rho_1) - \rho_2 }_1\\
\leq& \qty(1-\bar{\lambda})\qty(d_0 + \epsilon_1) + d_1\qty(\rho_1) + \epsilon_2\\
=& \qty(1-\bar{\lambda}) d_0  + d_1\qty(\rho_1) + \frac{\bar{\lambda}\epsilon'}{2}.
\end{align*}
Hence with $\rho_e = \rho_1$, the inequality \eqref{eq:d_Clambda_1} is satisfied. Notice that if $\mathcal{S}_{t}\qty(\rho) \in \mathscr{C}_{\ast}$, then $\rho \in \mathscr{C}_{\ast}$; see Proposition \ref{prp:C_astIsTwoSideInvariant}. Therefore, since $\rho_e \in \mathscr{C}_{\ast}$, there exists $\rho_c \in \mathscr{C}_\ast$ such that $\rho_e = \mathcal{S}_{t_n}\qty(\rho_c) \in \mathscr{C}_{\ast}$. Now with $\rho_m =  \bar{\lambda} \rho' + \qty(1-\bar{\lambda}) \rho_c$, the inequalities \eqref{eq:d_1} and \eqref{eq:d_Clambda_1} imply that
\begin{align}
0< \bar{\lambda} \epsilon' \leq \bar{\lambda} d\qty(\mathcal{S}_{t_n}\qty(\rho'),\mathscr{C}_\ast) \leq d\qty(\mathcal{S}_{t_n}\qty(\rho_m),\mathscr{C}_\ast) + \frac{\bar{\lambda}\epsilon'}{2} < \bar{\lambda} \epsilon'. \label{eq:d_Clambda_2}
\end{align}
Therefore we have arrived at a contradiction, which completes the proof; see also Figure \ref{fig:Global_Local_Asymptotic} for an illustration of the key ideas of this proof.\qed
\end{pf}
To this end, we can also highlight the following remarks:
\begin{rem}
	In contrast to the stability conditions given in \citep{Pan2014}, we do not require $V$ to be coercive, nor we demand that it commutes with the Hamiltonian of the system \citep{Pan2016}. We show in Theorem \ref{thm:LyapunovStability} that less restrictive conditions on both $V$ and $\mathcal{L}(V)$, \eqref{eq:V_rho_ast_0},\eqref{eq:AsymptoticallyStableCondition} are sufficient to guarantee the convergence of the density operator evolution to the set of invariant density operators.
\end{rem}
\begin{rem}
	We can use Theorem \ref{thm:LyapunovStability} to strengthen many results in the coherent control of quantum systems. In fact, the differential dissipative inequality given in \citep[Theorem 3.5]{James2010} and those which is given as an \gls{LMI} in \cite[Theorem 4.2]{James2008} explicitly imply global exponential and asymptotic stability conditions, provided that the storage function in \citep{James2010} and in \citep{James2008} have  global minima at the invariant density operator $\rho_\ast$. 
\end{rem}

\section{Stability of Quantum Systems II: An Invariance Principle}\label{sec:InvariantPrinciple}
In this section, we continue our discussion on the stability of the invariant density operator. This time we address the issue of analyzing the asymptotic stability condition of the set of invariant density operators when the generator of the Lyapunov observable fails to have a strict negative value in the neighborhood of the set of invariant density operators. This situation is similar to the Barbashin-Krasovskii-La Salle stability criteria for classical systems; see \citep{Krasovsky1963,khalil2002nonlinear}. In this section, we will use the weak topology to analyze the convergence of the density operator to the invariant density operators set. This choice is made since in the absence of the Lyapunov observable, we cannot use Corollary \ref{cor:WeakConverge} to infer convergence in $\norm{\cdot}_1$. We will rely in our analysis on some facts from functional analysis, to replace the compactness condition used in the invariance principle lemma by a more general condition in the context of a Banach space in the weak topology.
\\
Let us define a monotonically increasing sequence of times (MIST) $\qty{t_n}$ to be a sequence of times such that $0 \leq t_0 < t_1 < \cdots < t_n < \cdots$, and $t_n \rightarrow \infty$ as $n \rightarrow \infty$. In addition, we need the following definitions to establish the invariance principle.
\begin{defn}
	Let $\rho_t = \mathcal{S}_t\qty(\rho)$, where $\mathcal{S}_t$ is the predual semi-group \eqref{eq:SSemiGroup}. An element $\sigma \in \mathfrak{S}\qty(\mathscr{H})$ is said to be a positive limit of $\rho_t$, if there is an {MIST} $\qty{t_n}$, so that if $\rho_n \equiv \mathcal{S}_{t_n}\qty(\rho)$, $\rho_n \xrightharpoondown{w} \sigma$. We call the set $\mathscr{M}^+$ of all positive limits of $\rho_t$ the positive limit set of $\rho_t$.	
\end{defn}
\begin{defn}
	A set of density operators $\mathscr{M}$ is said to be invariant with respect to the {QDS} $\mathcal{T}_t$ \eqref{eq:TSemiGroup}, if for some $t \geq 0$, $\rho_t \in \mathscr{M}$ implies that $\rho_{t'} \in \mathscr{M}$ for any $t' \geq t$. Moreover, a set of density operators $\mathscr{M}'$ is said to be two-side invariant with respect to the {QDS} $\mathcal{T}_t$ if for some $t \geq 0$, $\rho_t \in \mathscr{M}'$, implies that $\rho_{t'} \in \mathscr{M}$ for any $t' \geq 0$.
\end{defn}
Notice the difference between an invariant set $\mathscr{M}'$, and a two-side invariant set $\mathscr{M}$. If $\mathscr{M}'$ is invariant then if $\rho_t \in \mathscr{M}'$ at time $t$, it will be contained in this set in the future. Moreover, a two-side invariant set $\mathscr{M}$ is an invariant set which also guarantees that $\rho_{t'}$ lies inside $\mathscr{M}$ for all the previous positive times before $t$. A set of invariant density operators $\mathscr{C}_\ast$ is by definition invariant, but it is also two-side invariant as shown by the following proposition.
\begin{prop}\label{prp:C_astIsTwoSideInvariant}
	The set $\mathscr{C}_\ast$ of invariant density operators of a semi-group $\mathcal{T}_t$ is two-side invariant.
\end{prop}
\begin{pf}
	Let $\mathcal{S}_t\qty(\rho) \in \mathscr{C}_\ast$, and let $t' \geq t$. By definition, for any $s \geq 0$ and any $A \in \mathfrak{B}\qty(\mathscr{H})$, $\trace\qty(\mathcal{S}_{t+s}\qty(\rho)A) =  \trace\qty(\mathcal{S}_{t}\qty(\rho)A)$. Therefore, for any $s\geq 0$,  $t'+s = t + (t'-t)+s$ where $(t'-t)+s \geq0$. Since $\trace\qty(\mathcal{S}_{t+(t'-t)}\qty(\rho)A) =  \trace\qty(\mathcal{S}_{t}\qty(\rho)A)$, we obtain for any $s \geq 0$, and any  $A \in \mathfrak{B}\qty(\mathscr{H})$ 
	\begin{align*}
	\trace\qty(\qty(\mathcal{S}_{t'+s}\qty(\rho) - \mathcal{S}_{t'}\qty(\rho))A) = 0.
	\end{align*}
	Now let $t' < t$. First, we claim that the semi-group $\mathcal{T}_t$ is a surjective mapping from $\mathfrak{B}\qty(\mathscr{H})$ onto $\mathfrak{B}\qty(\mathscr{H})$. To show this, for any $Y \in \mathfrak{B}\qty(\mathscr{H})$, there should be an $X \in \mathfrak{B}\qty(\mathscr{H})$ such that $\mathcal{T}_t\qty(X) = Y$. However,
	\begin{align*}
	Y =\mathcal{T}_t\qty(X) =& \rhoConditional{0]}{U_t^\dagger\qty(X\otimes \mathbb{1})U_t},\\
	Y \otimes \mathbb{1} =& U_t^\dagger\qty(X\otimes \mathbb{1})U_t,\\
	\rhoConditional{0]}{U_t\qty(Y\otimes \mathbb{1})U_t^\dagger} =& X.
	\end{align*}
	Therefore, since $\mathcal{S}_t\qty(\rho) \in \mathscr{C}_\ast$, we can write for any $s \geq 0$ and any $A \in \mathfrak{B}\qty(\mathscr{H})$
	\begin{align*}
	0 =& \trace\qty(\qty(\mathcal{S}_{t+s}\qty(\rho) - \mathcal{S}_{t}\qty(\rho))A) \\
	=& \trace\qty(\qty(\mathcal{S}_{(t-t')+t'+s}\qty(\rho) - \mathcal{S}_{(t-t')+t'}\qty(\rho))A)\\
	=& \trace\qty(\qty(\mathcal{S}_{t'+s}\qty(\rho) - \mathcal{S}_{t'}\qty(\rho))\mathcal{T}_{t-t'}\qty(A)).
	\end{align*}
	Due to the surjectivity of $\mathcal{T}_t$; i.e., $\mathcal{T}_t\qty(\mathfrak{B}\qty(\mathscr{H}))=\mathfrak{B}\qty(\mathscr{H})$, then $\trace\qty(\qty(\mathcal{S}_{t'+s}\qty(\rho) - \mathcal{S}_{t'}\qty(\rho))\mathcal{T}_{t-t'}\qty(A)) = 0$ for any $s \geq 0$ and any $A \in \mathfrak{B}\qty(\mathscr{H})$ is equivalent to $\trace\qty(\qty(\mathcal{S}_{t'+s}\qty(\rho) - \mathcal{S}_{t'}\qty(\rho))B) = 0$, for any $s \geq 0$ and any $B \in \mathfrak{B}\qty(\mathscr{H})$. This completes the proof.\qed	
\end{pf}
Before we begin to present the main result of this section, we first notice that in the classical invariance principle, the main asymptotic stability result was established via an invariance principle lemma, \citep[Lemma 4.1]{khalil2002nonlinear}. Roughly, it states that for any solution of a system of nonlinear differential equations which is bounded and belongs to some domain $\mathscr{D} \subseteq \mathbb{R}^n$, then its positive limit set $\mathscr{M}^+$ is non empty, compact and two-side invariant. Moreover, the solution approaches $\mathscr{M}^+$ as time goes to infinity. The proof of this lemma depends on the Bolzano-Weierstra\ss{} Theorem, which states that each bounded sequence in $\mathbb{R}^n$ has a convergent subsequence; see for example \citep[Thm 4.19]{Sutherland2009}. In what follow, we will discuss an equivalent condition of this theorem in the context of {QDS} in the weak topology. See Appendix for details on the weak topology used in this section. 
For an $A \in \mathfrak{B}\qty(\mathscr{H})$, and a $\rho \in \mathfrak{S}\qty(\mathscr{H})$ we define a semi-norm 
\begin{align}
p_A\qty(\rho) = \abs{\trace\qty(\rho A)}. \label{eq:pA}
\end{align}
It is obvious that a sequence of density operators $\qty{\rho_n}$ weakly converges to $\rho$, if and only if the limit $\lim\limits_{n\rightarrow \infty} p_A(\rho_n-\rho) = 0$ for all $A \in  \mathfrak{B}\qty(\mathscr{H})$. To examine the weak convergence of a sequence of density operators to a particular subset of $\mathfrak{S}\qty(\mathscr{H})$ we need to define the following. Let $\mathscr{C} \subset \mathfrak{S}\qty(\mathscr{H})$. We define $\rho_n \xrightharpoondown{w} \mathscr{C}$, if for all $A \in \mathfrak{B}\qty(\mathscr{H})$, $\lim\limits_{n\rightarrow \infty} \inf_{\sigma \in \mathscr{C}} p_A\qty(\rho_n-\sigma) = 0$.
\\
Now, we are ready to state the first result of this section. The following lemma is modification of the invariance principle lemma for quantum density operators.
\begin{lem}\label{lem:InvariantSetLemma}
	Suppose $\rho_t = \mathcal{S}_t\qty(\rho)$ where $\mathcal{S}_t$ is the predual semi-group \eqref{eq:SSemiGroup}. The positive limit set $\mathscr{M}^+$ of $\rho_t$ is nonempty, two-sided invariant, and weakly compact. Moreover, $\rho_t \xrightharpoondown{w} \mathscr{M}^+$.
\end{lem}
\begin{pf}
	Fix an initial density operator $\rho \in \mathfrak{S}\qty(\mathscr{H})$. Since $\mathfrak{S}\qty(\mathscr{H})$ is a weakly compact subset of $\mathfrak{I}_1\qty(\mathscr{H})$, by Lemma \ref{lem:EberleinSmulian} in the Appendix, every sequence in $\mathfrak{S}\qty(\mathscr{H})$ has a subsequence converging weakly to an element of $\mathfrak{S}\qty(\mathscr{H})$. Therefore, for any {MIST} $\qty{t_n}$, by defining $\rho_n \equiv \mathcal{S}_{t_n}\qty(\rho)$, then $\qty{\rho_n}$ is a sequence in $\mathfrak{S}\qty(\mathscr{H})$, which has a weakly convergent subsequence $\qty{\rho_{n_k}}$ converging to an element $\rho_c \in \mathfrak{S}\qty(\mathscr{H})$; i.e., $\rho_{n_k} \xrightharpoondown{w} \rho_c$. Therefore, the positive limit set $\mathscr{M}^+$ is nonempty.
	\\
	To show that $\mathscr{M}^+$ is a two-side invariant set, 
	let us first assume that $t' \geq t$. Suppose $\mathcal{S}_t\qty(\rho) \in \mathscr{M}^+$, then there exists an {MIST} $\qty{t_n}$ such that $\mathcal{S}_{t_n}\xrightharpoondown{w} \mathcal{S}_t\qty(\rho)$. Therefore, by the semi-group property, we may write, $S_{t'}\qty(\rho) = \mathcal{S}_{(t'-t)}\qty(\mathcal{S}_t\qty(\rho)) = w-\lim\limits_{n \rightarrow \infty} \mathcal{S}_{(t'-t)}\qty(\mathcal{S}_{t_n} \qty(\rho))$. Therefore there exists an {MIST} $\qty{\tau_n}$, where $\tau_n = (t'-t)+t_n$ such that  $w-\lim\limits_{n \rightarrow \infty} \mathcal{S}_{\tau_n} \qty(\rho) = \mathcal{S}_{t'}\qty(\rho)$. Therefore $\mathcal{S}_{t'}\qty(\rho) \in \mathscr{M}^+$.\\
	Now assume that $t' < t$, $t'\geq0$. Suppose $\mathcal{S}_t\qty(\rho) = w-\lim\limits_{n \rightarrow \infty} \mathcal{S}_{t_n} \qty(\rho)$. Then $\mathcal{S}_{t'}\qty(\rho) = \mathcal{S}_{t'-t}\qty(\mathcal{S}_{t}\qty(\rho)) =  w-\lim\limits_{n \rightarrow \infty} \mathcal{S}_{t'-t}\qty(\mathcal{S}_{t_n} \qty(\rho))$. As $t_n \rightarrow \infty$, we always have a subsequence $\qty{t_m}$ of $\qty{t_n}$ such that $(t'-t)+t_m \geq 0$. Therefore, starting at a sufficiently large $n$ as the first element of the subsequence $\qty{t_m}$, there exists an {MIST} $\qty{\tau_m}$, where $\tau_m = (t'-t)+t_m$ such that  $w-\lim\limits_{m \rightarrow \infty} \mathcal{S}_{\tau_m} \qty(\rho) = \mathcal{S}_{t'}\qty(\rho)$. It follows that $\mathcal{S}_{t'}\qty(\rho) \in \mathscr{M}^+$. Hence, we conclude that $\mathscr{M}^+$ is a two-side invariant set.
	\\
	To show that $\mathscr{M}^+$ is a weakly compact set, it is sufficient to show that $\mathscr{M}^+$ is a weakly closed subset of $\mathfrak{S}\qty(\mathscr{H})$.  Suppose that there exists a net of density operators $\qty{\rho_{i,c}} \in \mathscr{M}^+$, weakly converging to a density operator $\rho_c \in \mathfrak{S}\qty(\mathscr{H})$, $\rho_{i,c}\xrightharpoondown{w}\rho_c$. Let $\qty{t_{i,n}}$ be an {MIST} such that
	$\rho_{i,n} = \mathcal{S}_{t_{i,n}}\qty(\rho)$ weakly converges to $\rho_{i,c}$. Then we have to show that there exists an {MIST} $\qty{\tau_i}$ such that $\rho_{i} = \mathcal{S}_{\tau_{i}\qty(\rho)}$ weakly converges to $\rho_c$. Suppose that this statement is false. Then for any {MIST} $\qty{t_n}$ there exists a subsequence $\qty{t_m}$ of $\qty{t_n}$ such that for any $m$, $p_A\qty(\mathcal{S}_{t_m}\qty(\rho) - \rho_c )\geq \epsilon$ for $\epsilon >0$ and $A \in \mathfrak{B}\qty(\mathscr{H})$. However, for the subsequence $\qty{t_m}$, there exists an element $\rho_{i,c} \in \mathscr{M}^+$  such that $\mathcal{S}_{t_m}\qty(\rho) \xrightharpoondown{w} \rho_{i,c}$. We can write for any $m$, 
$
	{0 < \epsilon} \leq p_A\qty(\mathcal{S}_{t_{m}}\qty(\rho) - \rho_c ) \leq p_A\qty(\mathcal{S}_{t_{m}}\qty(\rho) - \rho_{i,c} ) + p_A\qty(\rho_{i,c} - \rho_c ).$
	Since the net $\rho_{i,c} \xrightharpoondown{w}\rho_c$ and $\mathcal{S}_{t_m}\qty(\rho) \xrightharpoondown{w} \rho_{i,c}$ the right hand side of the above equation will go to zero. Hence we establish a contradiction.
	\\
	Finally, we want to show that $\rho_t \xrightharpoondown{w} \mathscr{M}^+$. Again, we will prove this by a contradiction. Suppose this statement is false. Then there exists an {MIST} $\qty{t_n}$ such that $\inf_{\sigma \in \mathscr{M}^+} p_A \qty(S_{t_n}\qty(\rho),\sigma) \geq \epsilon$ for any $n$ and some $\epsilon >0$, and an $A \in \mathfrak{B}\qty(\mathscr{H})$. But  $\mathfrak{S}\qty(\mathscr{H}) $ is weakly compact. Therefore, $S_{t_n}\qty(\rho)$ contains a weakly convergent subsequence $S_{t_m}\qty(\rho)$. Hence, there exists $\rho_c \in \mathfrak{S}\qty(\mathscr{H})$ such that $w-\lim\limits_{n\rightarrow \infty} S_{t_m}\qty(\rho) = \rho_c$. Therefore, $\rho_c$ must be an element of $\mathscr{M}^+$, but at the same time, for some $\epsilon >0$ and $A \in \mathfrak{B}\qty(\mathscr{H})$, we have
$
	\epsilon < \lim\limits_{m\rightarrow \infty} \inf_{\sigma \in \mathscr{M}^+} p_A \qty(S_{t_m}\qty(\rho),\sigma) \leq  \lim\limits_{m\rightarrow \infty} p_A \qty(S_{t_m}\qty(\rho),\rho_c)  = 0,
	$
	which is a contradiction. Therefore 
	\begin{align*}
	\lim\limits_{n\rightarrow \infty} \inf_{\sigma \in \mathscr{M}^+} p_A \left(S_{t_n}\qty(\rho),\sigma\right)=0,
	\end{align*}
	for $A \in \mathfrak{B}\qty(\mathscr{H})$ which completes the proof.\qed
\end{pf}
\begin{rem}
	If the positive limit set $\mathscr{M}^+$ consists only a unique element of $\mathfrak{S}\qty(\mathscr{H})$, then this element is an invariant density operator and it is weakly attractive globally; see also \citep[Thm 1]{Schirmer2010}. 
\end{rem}
Using Lemma \ref{lem:InvariantSetLemma}, we are ready to state a quantum version of Barbashin-Krasovskii-La Salle's Theorem for invariant density operators.
\begin{thm}\label{thm:LaSalle}
	Let $V \in \mathfrak{B}\qty(\mathscr{H})$ be a self-adjoint operator such that for some neighborhood $\mathscr{N}$, $\trace\qty(\mathcal{L}\qty(V)\rho) \leq 0, \; {\forall \rho \in \mathscr{N}}$. Also, let $\mathscr{C}$ be a weakly compact subset of $\mathscr{N}$ that is invariant. Define a set $\mathscr{E}$ as follows:
	\begin{align}
	\mathscr{E} = \qty{\rho \in \mathscr{C}: \trace\qty(\mathcal{L}\qty(V)\rho) = 0}.
	\end{align}
	In addition, let $\mathscr{M}$ be the largest two-side invariant set in $\mathscr{E}$. Then for any $\rho \in \mathscr{C}$, $\mathcal{S}_t\qty(\rho)\xrightharpoondown{w}\mathscr{M}$.
\end{thm}
\begin{pf}
	Let $\rho_t = \mathcal{S}_t\qty(\rho)$, where initially $\rho \in \mathscr{C}$. Since $\mathscr{C}$ is invariant, any $\rho_t$ starting in it, will remain inside it in the future. Therefore, the positive limit set $\mathscr{M}^+$ of $\rho_t$ lies inside $\mathscr{C}$ since $\mathscr{C}$ is weakly compact. Observe that $\trace\qty(\mathcal{L}\qty(V)\rho) \leq 0$ on $\mathscr{N}$. Therefore, $\trace\qty(V \rho_t)$ is non-increasing; that is for any $s<t$, $\trace\qty(V \rho_t) \leq \trace\qty(V \rho_s)$. Since $\trace\qty(V \rho_t): \mathfrak{I}_1\qty(\mathscr{H}) \rightarrow \mathbb{C}$ is weakly continuous and $\mathscr{C}$ is weakly compact, then $\trace\qty(V \rho_t)$ is bounded from below on $\mathscr{C}$. Therefore, there is a $v \in \mathbb{R}$ such that $\lim\limits_{t\rightarrow \infty} \trace\qty(V \rho_t) = v$. For any $\rho_c \in \mathscr{M}^+$, there is an MIST such that $w-\lim\limits_{n\rightarrow \infty} \mathcal{S}_{t_n}\qty(\rho) = \rho_c$. Since $\trace\qty(V \rho_t)$ is weakly continuous, then $\lim\limits_{n\rightarrow \infty} \Paverage{V_{t_n}} = \lim\limits_{n\rightarrow \infty} \trace\qty(V \mathcal{S}_{t_n}\qty(\rho)) = \trace\qty(V \rho_c) = v$. Therefore, on $\mathscr{M}^+$, $\trace\qty(V \rho_t) = v$. By Lemma \ref{lem:InvariantSetLemma}, $\mathscr{M}^+$ is a two-side invariant set. Therefore, $\trace\qty(\mathcal{L}\qty(V)\rho) = 0$ on $\mathscr{M}^+$. Consequently, we conclude that
	\begin{align*}
	\mathscr{M}^+ \subset \mathscr{M} \subset \mathscr{E} \subset \mathscr{C}.
	\end{align*}
	Since $\mathscr{C}$ is by definition weakly compact, then Lemma \ref{lem:InvariantSetLemma} implies that $\rho_t \xrightharpoondown{w} \mathscr{M}^+$. Therefore, $\rho_t \xrightharpoondown{w} \mathscr{M}$. \qed		
\end{pf}
The following corollary of Theorem \ref{thm:LaSalle} shows that in the case that the generator of a Lyapunov candidate operator $V$ fails to have a strictly negative value, as long as the invariant density operator $\mathscr{C}_\ast$ is the largest two-sided invariant set having $ \trace\qty(\mathcal{L}\qty(V)\rho) = 0$, then we can still infer weak asymptotic stability.
\begin{cor}\label{cor:LaSalleAsymptoticStability}
	Suppose the set of invariant density operators of the {QDS} $\mathcal{T}_t$ \eqref{eq:TSemiGroup} is given by $\mathscr{C}_\ast$.  Let $V \in \mathfrak{B}\qty(\mathscr{H})$ be a self-adjoint operator such that for some weakly compact invariant set $\mathscr{C}$ for which $\mathscr{C}_\ast \subset \mathscr{C}$,  $\trace\qty(\mathcal{L}\qty(V)\rho) \leq 0, \; {\forall \rho \in \mathscr{C}}$. Define a set $\mathscr{E}$ as follows:
	\begin{align}
	\mathscr{E} = \qty{\rho \in \mathscr{C}: \trace\qty(\mathcal{L}\qty(V)\rho) = 0}.
	\end{align}
	Suppose that the largest two-side invariant set in $\mathscr{E}$ is $\mathscr{C}_\ast$.  Then for any $\rho \in \mathscr{C}$, $\rho_t \xrightharpoondown{w} \mathscr{C}_\ast$.
\end{cor}

\section{Examples}
To give an indication of the applications of the Lyapunov stability and the invariance principle conditions we have derived, we consider the following examples. 
\subsection{Displaced Linear Quantum System}\label{sec:exm_LinearQuantumSystem}
	Consider a linear quantum system $\mathcal{P}$, with $\mathbb{H} =  (a - \alpha \mathbb{1})^\dagger (a - \alpha \mathbb{1}), \alpha \in \mathbb{C}$, $\mathbb{L} = \sqrt{\kappa} \qty(a - \alpha\mathbb{1})$, and $\mathbb{S} = \mathbb{1}$. Evaluating the steady state of \eqref{eq:MasterEquation} we know that the invariant density operator is a coherent density operator with amplitude $\alpha$; i.e., $\rho_{\ast} = \ket{\alpha}\bra{\alpha}$. Now, choose the Lyapunov observable $V = \mathbb{H}$. Straightforward calculation of $\mathcal{L}(V)$ \citep{gough2009series} gives, 
	\begin{equation*}
	\mathcal{L}\qty(V) = -\kappa (N -\frac{1}{2}\qty(\alpha a^\dagger + \alpha^\ast a ) + \abs{\alpha}^2 \mathbb{1}). 
	\end{equation*}
	Notice that $\trace\qty(\mathcal{L}\qty(V)\rho)< 0$ for all density operators other than $\rho_{\ast} = \ket{\alpha}\bra{\alpha}$. To verify this inequality, it is sufficient to take $\rho_t = \ket{\beta}\bra{\beta}$ with $\beta \neq \alpha$. 
	Therefore, for $\rho_t = \ket{\beta}\bra{\beta}$, we obtain $\trace \qty(\mathcal{L}\qty(V) \rho_t) = \abs{\beta}^2 + \abs{\alpha}^2 -\qty(\alpha^\ast \beta+\beta^\ast \alpha)$. Therefore, $\trace\qty(\mathcal{L}\qty(V)\rho_t) = -\kappa \trace \qty(V \rho_t) + \kappa/2 \qty(\alpha^\ast \beta+\beta^\ast \alpha) <0$. Hence, 
Theorem \ref{thm:LyapunovStability} indicates that the invariant density operator is exponentially stable. %

\subsection{Displaced Linear Quantum System, with Alternative Lyapunov Operators}\label{sec:exm_LinearQuantumSystemWithOtherV}
	Consider again the linear quantum system $\mathcal{P}$ in Example \ref{sec:exm_LinearQuantumSystem}, where this time we set $\alpha = 0$ and $\kappa = 1$. From the previous example, the invariant density operator  is the vacuum state; i.e., $\rho_\ast = \ket{0}\bra{0}$ . Suppose instead of selecting the Lyapunov operator as in the previous example, we select $V = N^2 - N$. We then obtain
	\begin{equation*}
	{\mathcal{L}\qty(V) =  -\qty(2N\qty(N - \mathbb{1}))}. 
	\end{equation*}
	Using number state $\ket{n}$, we obtain $\trace\qty(\mathcal{L}\qty(V)\rho) =-\qty( 2n\qty(n-1) ) \leq 0$. For any $\ket{n}$, the generator of the Lyapunov observable will be zero when $n = 0$ and $n = 1$. We further notice that if we choose the neighborhood of $\rho_\ast$ to be the whole space of density operators; i.e., $\mathscr{N} = \mathfrak{S}\qty(\mathscr{H})$, the strict minima condition in \eqref{eq:V_rho_ast_0} is not satisfied since when $\rho_t \in \mathscr{R} \equiv \qty{\lambda \rho_\ast + \qty(1-\lambda) \ket{1}\bra{1}: \lambda\in (0,1]}$, $\trace\qty(V \rho) = 0$ and $\trace\qty(\mathcal{L}\qty(V)\rho) = 0$. Reducing the neighborhood of $\rho_\ast$ so that the strict minima condition in \eqref{eq:V_rho_ast_0} is satisfied is also impossible, since for any $\epsilon >0$,  selecting $\mathscr{N} = \mathscr{B}_{\epsilon}\qty(\rho_\ast)$ will have intersection with $\mathscr{R}$. Therefore, with this Lyapunov candidate we cannot use Theorem \ref{thm:LyapunovStability} to analyze the stability of this invariant density operator. However, we can use our invariance principle to infer the global weak asymptotic stable property of $\rho_\ast$. Suppose we choose $\mathscr{C} = \mathfrak{S}\qty(\mathscr{H})$ in Corollary \ref{cor:LaSalleAsymptoticStability}. According to this corollary, since the largest two-side invariant set on the set $\mathscr{E} = \qty{\rho \in \mathscr{C}: \trace\qty(\mathcal{L}\qty(V)\rho) = 0}$, is given by $\rho_\ast$. Therefore, $\rho_\ast$ is globally weakly asymptotically stable.

\subsection{Nonlinear Quantum System With Non Unique Equilibrium Points}\label{sec:exm_NonLinearQuantumSystem}
\begin{figure}[h]
		\centering
		\includegraphics[width=0.5\textwidth]{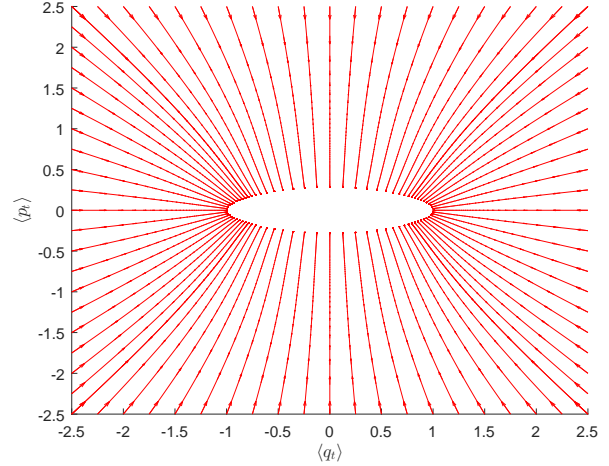}
		\caption{Illustration of trajectories of the quantum system in Example \ref{sec:exm_NonLinearQuantumSystem} which are simulated using the corresponding master equation. Each curve corresponds to a different initial coherent density operator. 
		} 
		\label{fig:SMESimulation}
	\end{figure}

\begin{figure}[h]
		\centering
		\includegraphics[width=0.5\textwidth]{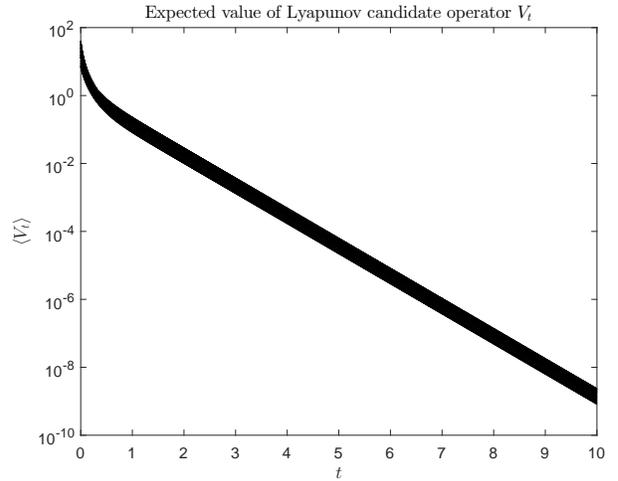}
		\caption[Lyapunov operator expected value of the quantum system in Example \ref{sec:exm_NonLinearQuantumSystem}.]{Lyapunov operator expected value of the quantum system in Example \ref{sec:exm_NonLinearQuantumSystem}.} 
		\label{fig:Example3Lyapunov}
	\end{figure}

	In this example, we consider a quantum system where the set of invariant states contain more than one state. Experimental stabilization of such a system has been reported recently in \citep{Leghtas2015}. 
	Consider a nonlinear quantum system with zero Hamiltonian and a coupling operator $\mathbb{L} = (a^2 - \alpha^2\mathbb{1}) $,  where $\alpha$ is a complex constant; see also \citep{Mirrahimi2014}. 
	To find the invariant density operators of this quantum system we need to find the eigenvectors of $a^2$. Without loss of generality, let $\ket{z}$ be one of the eigenvectors of $a^2$, such that $a^2\ket{z} = \alpha^2 \ket{z}$. Expanding $\ket{z}$ in the number state orthogonal basis, we can write, $a^2 \ket{z} = \sum_{n=0}^{\infty}  a^2 c_n \ket{n}$. Therefore, we find that, $\alpha^2 c_{n-2} = c_n \sqrt{n\qty(n-1)}$. By mathematical induction, we have for $n$ even, $c_n = c_0 \alpha^n/\sqrt{n!}$, and for $n$ odd, $c_n = c_1 \alpha^n/\sqrt{n!}$. Therefore, we can write the eigenvector of $a^2$ as,
	\begin{align*}
	\ket{z} =& c_0 \sum_{n=0, n \text{even}}^{\infty} \dfrac{\alpha^n}{\sqrt{n!}} \ket{n} + c_1 \sum_{n=1, n \text{odd}}^{\infty} \dfrac{\alpha^n}{\sqrt{n!}} \ket{n}.
	\end{align*}
	By observing that a coherent vector with magnitude $\alpha$, is given by 
	\begin{align*}
	\ket{\alpha} = \exp(-\dfrac{\abs{\alpha}^2}{2})\sum_{n=0}^{\infty} \dfrac{\alpha^{n}}{\sqrt{n!}} \ket{n},
	\end{align*}
	we can write
	\begin{align*}
	\ket{z} =& c_0 \exp(\frac{\abs{\alpha}^2}{2}) \frac{\qty(\ket{\alpha} + \ket{-\alpha})}{2}\\
	 		+& c_1 \exp(\frac{\abs{\alpha}^2}{2}) \frac{\qty(\ket{\alpha} - \ket{-\alpha})}{2}.
	\end{align*}
	Normalization of $\ket{z}$ shows that $c_0$ and $c_1$ satisfy an elliptic equation,
	\begin{align*}
	\abs{c_0}^2 \cosh(\abs{\alpha}^2)+ \abs{c_1}^2 \sinh(\abs{\alpha}^2) = 1.
	\end{align*} 
	Therefore, we can construct	the set of solutions of $\mathbb{L} \ket{z} = 0 \ket{z}$ as the following set:
	\begin{equation}
	\mathscr{Z}_\ast = \qty{\ket{z}\in \mathscr{H}: \ket{z} = C_0\frac{\qty(\ket{\alpha} + \ket{-\alpha})}{2} + C_1\frac{\qty(\ket{\alpha} - \ket{-\alpha})}{2}}, 
	\end{equation}
	where $C_0 = c_0 \exp(\frac{\abs{\alpha}^2}{2})$ and $C_1 = c_1 \exp(\frac{\abs{\alpha}^2}{2})$. The set of invariant density operators of this quantum system is a convex set $\mathscr{C}_\ast$ which is given by
	\begin{equation}
	\mathscr{C}_\ast = \qty{\sum_i \lambda_i \ket{\beta_i}\bra{\beta_i}:   \ket{\beta_i} \in \mathscr{Z}_\ast}, 
	\end{equation}
	where $\lambda_i \geq 0, \sum_i \lambda_i = 1$. Suppose we select a Lyapunov candidate $V = \mathbb{L}^\dagger \mathbb{L}$. One can verify that $\trace\qty(\rho V) = 0$ for all $\rho$ belonging to $\mathscr{C}_\ast$, and has a positive value outside of this set.
	\\
	Straightforward calculation of the quantum Markovian generator of $V$
	 gives us the following
	\begin{equation}
	\mathcal{L}\qty(V) = - \qty(4 \mathbb{L}^\dagger N \mathbb{L} + 2 V).
	\end{equation}
	Outside the set $\mathscr{C}_\ast$, the generator $\mathcal{L}\qty(V)$ has a negative value.  Therefore, Theorem \ref{thm:LyapunovStability} implies that the set $\mathscr{C}_\ast$ is globally exponentially stable. Convergence analysis of quantum system in this example has also been analyzed in \citep{Azouit2015} using different metric, assuming that initially the system's density operator is spanned by a finite collection of number states $\ket{n}$.\\
	Figure \ref{fig:SMESimulation} illustrates the phase-space of the system corresponding to various initial density operators.  It reveals that each distinct trajectory converges to a different invariant density operator, all belonging to the set of invariant density operators $\mathscr{C}_\ast$. 
	We can also calculate the quantum Markovian generator of $q_t$ and $p_t$ as below
	\begin{align}
	\mathcal{G}\qty(\mathbf{x}_t) =& -\frac{1}{4}\mqty[2q_t^3+p_t^2q_t+q_tp_t^2-4q_t\\2p_t^3+q_t^2p_t+p_tq_t^2-4p_t] \nonumber \\
	+& \mqty[\Re\qty(\alpha^2) && \Im\qty(\alpha^2)\\ -\Im\qty(\alpha^2) && -\Re\qty(\alpha^2)]\mqty[q_t\\p_t] \label{eq:NonlinearDynamicsExample.5.3}
	\end{align}
	Moreover, Figure \ref{fig:Example3Lyapunov} shows the Lyapunov operator expected values. It also explains that although each trajectory converges to a distinct invariant density operator, their Lyapunov expected values do converge to zero. Such a convergence indicates once more that the invariant density operators are in the set $\mathscr{C}_\ast$.

\section{Conclusion}
In this article, we have proposed a Lyapunov stability notion for open quantum systems, which enables the analysis of the convergence of the system's density operator in $\norm{ \cdot }_1$. Under a uniform continuity assumption, this stability notion is stronger compared to the weak convergence and finite moment convergence that have been considered for quantum systems previously. 
\\
We have shown how to analyze the stability of this set via a candidate Lyapunov operator. We have also shown that in uniformly continuous QDS, local asymptotic stability is equivalent to global asymptotic stability. Lastly, we have introduced an analog of the Barbashin-Krasovskii-La Salle Theorem on the dynamics of quantum systems, which brings the possibility to infer asymptotic stability using a Lyapunov candidate operator, even if its generator fails to have a strictly negative value in a neighborhood of the invariant density operators.

\section{Acknowledgments}
The authors acknowledge discussions with Dr. Hendra Nurdin of UNSW.

\section*{Appendix}
\subsection*{Weak Topology}\label{sec:WeakTopology}
The topology of a set determines how each member of this set relates to the others. It determines whether a set is open, or closed. In this subsection, we will be interested in a type of topology known as the \emph{weak} topology. First, we observe that the collection of all bounded linear functionals on a vector space $\mathscr{X}$, $\mathfrak{B}\qty(\mathscr{X},\mathbb{C})$ is also a vector space. For this collection, we write $\mathscr{X}^\ast = \mathfrak{B}\qty(\mathscr{X},\mathbb{C})$. We also say that $\mathscr{X}^\ast$ is the \emph{dual} of $\mathscr{X}$.   We use the same asterisk symbol $^\ast$ here to denote the dual $\mathscr{X}$. We notice that by using the Riez representation theorem \citep[Thm I.3.4]{Conway1985}, for the case of a Hilbert space, any bounded linear functional $\phi(\ket{\psi})$ always corresponds to a vector $\bra{\phi}$, such that $\phi(\ket{\psi}) = \expval{\phi,\psi}$. It is well known that if $\mathscr{X}$ is a normed space (not necessarily complete) then $\mathscr{X}^\ast$ is a Banach space \citep[Prop III.5.4]{Conway1985}. The \emph{weak} topology on $\mathscr{X}$ or $\sigma\qty(\mathscr{X},\mathscr{X}^\ast)$ is defined as the weakest topology on $\mathscr{X}$ that makes all elements in $\mathscr{X}^\ast$ continuous. We observe that for $x \in \mathscr{X}$, and $l \in \mathscr{X}^\ast$, then the element $x$ itself can be considered as a linear functional in $\mathscr{X}^\ast$ by defining $x(l) \equiv l(x)$. Hence, we also have a \emph{weak}-$^\ast$ topology or $\sigma\qty(\mathscr{X}^\ast,\mathscr{X})$, the topology on  $\mathscr{X}^\ast$ such that the functional $l(x)$ is continuous, for all $x \in \mathscr{X}$, and $l \in \mathscr{X}^\ast$.
\\
The convergence of a net $\qty{x_i}$ to an element $x_o$ in the weak topology can be characterized by the convergence of $l(x_i) \rightarrow l(x_o)$ for all $l \in \mathscr{X}^\ast$. Every element $Y$ in the set of bounded operators  $\mathfrak{B}\qty(\mathscr{H})$ defines a linear functional $\omega_Y$ on the set of trace bounded operators $\mathfrak{I}_1\qty(\mathscr{H})$ where $\omega_Y\qty(X) = \trace\qty(XY)$. 
Since $\mathfrak{B}\qty(\mathscr{H})$ is the dual of $\mathfrak{I}_1\qty(\mathscr{H})$, (or similarly $\mathfrak{I}_1\qty(\mathscr{H})$ is the \emph{predual} of $\mathfrak{B}\qty(\mathscr{H})$) \citep[Thm III.19.2]{conway2000course}, then letting $\mathscr{X} = \mathfrak{I}_1\qty(\mathscr{H})$, and $\mathscr{X}^\ast = \mathfrak{B}\qty(\mathscr{H})$, then the topology on $\mathfrak{B}\qty(\mathscr{H})$ defined by $\mathfrak{I}_1\qty(\mathscr{H})$, $\sigma\qty(\mathfrak{B}\qty(\mathscr{H}),\mathfrak{I}_1\qty(\mathscr{H}))$ is a weak $^\ast$ topology on $\mathscr{X}^\ast = \mathfrak{B}\qty(\mathscr{H})$. This topology has also been referred to with many different names, such as the \emph{ultraweak}, or \emph{$\sigma$ weak}, or \emph{weak $^\ast$ operator 
} or \emph{normal} topologies; see \citep[\textsection VI.6 p. 213]{Reed1972}, \citep[\textsection III.20]{conway2000course}, \citep[Appendix 1 p. 214]{Meyer1993}. In this article, we will refer to  $\sigma\qty(\mathfrak{B}\qty(\mathscr{H}),\mathfrak{I}_1\qty(\mathscr{H}))$ as the normal topology.
We will also use the topology $\sigma\qty(\mathfrak{I}_1\qty(\mathscr{H}),\mathfrak{B}\qty(\mathscr{H}))$. This topology is the weak topology on $\mathfrak{I}_1\qty(\mathscr{H})$. Therefore, to avoid an unnecessary notational complexity, we will use the \emph{weak} topology notion for $\sigma\qty(\mathfrak{I}_1\qty(\mathscr{H}),\mathfrak{B}\qty(\mathscr{H}))$. 
\begin{defn}\label{def:NormalConvergence}
	Let $\qty{X_n}$ be a sequence in $\mathfrak{B}\qty(\mathscr{H})$. We say $\qty{X_n}$ converges \emph{normally} to $X$ if for every $\rho$ in $\mathfrak{I}_1\qty(\mathscr{H})$, 
	\begin{dmath}
		\lim\limits_{i\rightarrow \infty} \trace\qty(X_i \rho) = \trace\qty(X \rho). \label{eq:NormalConvergence}
	\end{dmath}
\end{defn}
\begin{defn}\citep{Fagnola2003}\label{def:weakConvergence}
	A sequence of density operators $\qty{\rho_i}$ is said to converge \emph{weakly} to $\rho \in \mathfrak{S}\qty(\mathscr{H})$ if for all $A \in \mathfrak{B}\qty(\mathscr{H})$
	\begin{dmath}
		\lim\limits_{i \rightarrow \infty} \trace\qty(\rho_i A) = \trace\qty(\rho A). \label{eq:weakConvergence}
	\end{dmath}
\end{defn}
We write the limit of a weakly-convergent sequence $\qty{\rho_i}$ to $\rho$ as $w-\lim\limits_{i\rightarrow \infty} \rho_i = \rho$, or $\rho_i \xrightharpoondown{w} \rho$. 
These various convergences also apply directly without modification to generalized sequences (nets), which are not  necessary unless if we consider non-separable Hilbert spaces. This rarely happens in mathematical physics, \citep{davies1976quantum}. In the case that $\mathscr{X}$ is a normed space (not necessarily complete), then the norm and the weak topologies on $\mathscr{X}$ are equivalent if and only if $\mathscr{X}$ is finite-dimensional, \citep[Prop 2.5.13]{Megginson1998}.
For a Banach space $\mathscr{X}$, let us denote the closed unit ball in $\mathscr{X}$ by $\mathscr{B}_{\mathscr{X}} \equiv \qty{x\in\mathscr{X}: \norm{x}\leq 1}$. Generally, for an infinite dimensional Banach space $\mathscr{X}$, the unit ball $\mathscr{B}_{\mathscr{X}}$ is not compact in the topology induced by its norm. However, there are some results in analysis that imply the compactness of the unit ball in the weak topology. One such result is known as the Banach-Alaoglu Theorem. We will use this theorem as well as the Eberlein-\v{S}mulian Theorem to establish a positive invariance principle for {QDS}.
Using the Banach-Alaoglu Theorem \citep[Thm III.3.1 ]{Conway1985}, we obtain the following weak compactness result for $\mathfrak{S}\qty(\mathscr{H})$. 
\begin{cor}
	The set of density operators $\mathfrak{S}\qty(\mathscr{H})$ is weakly compact.
\end{cor}
\begin{pf}
	Let $\mathscr{X}$ be $\mathfrak{I}_1\qty(\mathscr{H})$. Then according to Banach-Alaoglu Theorem, the ball  $\mathscr{B}_{\mathscr{X}^{\ast\ast}}$ is weakly compact in the $\sigma\qty(\mathscr{X}^{\ast\ast},\mathscr{X}^\ast)$ topology. In this case, $\mathscr{X}^{\ast\ast} = \mathfrak{B}\qty(\mathscr{H})^\ast$. For each $x \in \mathscr{X}$, we can define a linear functional $\omega_x\qty(\cdot) = \trace\qty(x \cdot)$. Therefore, there exists a mapping $\iota : x \rightarrow \omega_x$ which is a linear isometry \citep[II.A.10]{Wojtaszczyk1991}. We can identify $\sigma\qty(\mathscr{X},\mathscr{X}^\ast)$ as a restriction of $\sigma\qty(\mathscr{X}^{\ast\ast},\mathscr{X}^\ast)$ to $\iota\qty(\mathscr{X})$; i.e., $\iota$ is a homeomorphism from $\mathscr{X}$ in the $\sigma\qty(\mathscr{X},\mathscr{X}^\ast)$ topology onto $\iota\qty(\mathscr{X})$ in the $\sigma\qty(\mathscr{X}^{\ast\ast},\mathscr{X}^\ast)$ topology \citep[II.A.10]{Wojtaszczyk1991}. Observe that $\mathfrak{S}\qty(\mathscr{H}) \subset \mathscr{X}$, and therefore, $\iota\qty(\mathfrak{S}\qty(\mathscr{H})) \subset \mathscr{B}_{\mathscr{X}^{\ast\ast}}$, since for any $\rho \in  \mathfrak{S}\qty(\mathscr{H})$, $\norm{\omega_\rho}\leq 1$. According to Proposition \ref{prp:SHisCompleteSubspaceOfI1}, $\mathfrak{S}\qty(\mathscr{H})$ is closed in the Banach space $\qty(\mathscr{X},\norm{\cdot}_1)$. Since $\mathfrak{S}\qty(\mathscr{H})$ is convex, $\mathfrak{S}\qty(\mathscr{H})$ is also weakly closed. This implies that $\iota\qty(\mathfrak{S}\qty(\mathscr{H}))$ is closed in $\sigma\qty(\mathscr{X}^{\ast\ast},\mathscr{X}^\ast)$, which implies that $\iota\qty(\mathfrak{S}\qty(\mathscr{H}))$ is compact in the $\sigma\qty(\mathscr{X}^{\ast\ast},\mathscr{X}^\ast)$ topology, so is $\mathfrak{S}\qty(\mathscr{H})$ in the $\sigma\qty(\mathscr{X},\mathscr{X}^\ast)$ topology.\qed
\end{pf}
Let us recall some notions on compactness \citep[\textsection 1.6]{Albiac2006}. If $\mathscr{X}$ is a topological space and $\mathscr{A} \subseteq \mathscr{X}$, $\mathscr{A}$ is said to be \emph{sequentially compact} if every sequence in $\mathscr{A}$ has a subsequence that converges to a point in $\mathscr{A}$ (to a point in $\mathscr{X}$, respectively). 
While sequential compactness is equivalent to compactness in a metric space, since the weak topology of an infinite dimensional vector space is not metrizable, the equivalence is no longer preserved. However, the following Lemma asserts that a Banach space in the weak topology behaves like a metric space, although it does not have to be metrizable. This result also shows that every sequence in $\mathfrak{S}\qty(\mathscr{H})$ has a weakly convergent subsequence. 
\begin{lem}[Eberlein-\v{S}mulian Theorem]\citep[Thm 1.6.3]{Albiac2006} \label{lem:EberleinSmulian}
	Let $\mathscr{X}$ be a Banach space. and let $A \subseteq \mathscr{X}$. The following conditions are then equivalent:
	\begin{enumerate}
		\item $A$ is weakly compact.
		\item $A$ is weakly sequentially compact.
	\end{enumerate}	
\end{lem}

\bibliographystyle{elsarticle-harv}        				  %
\bibliography{ReferenceAbbrvBibLatex}

\begin{thebibliography}{75}
\expandafter\ifx\csname natexlab\endcsname\relax\def\natexlab#1{#1}\fi
\providecommand{\url}[1]{\texttt{#1}}
\providecommand{\href}[2]{#2}
\providecommand{\path}[1]{#1}
\providecommand{\DOIprefix}{doi:}
\providecommand{\ArXivprefix}{arXiv:}
\providecommand{\URLprefix}{URL: }
\providecommand{\Pubmedprefix}{pmid:}
\providecommand{\doi}[1]{\href{http://dx.doi.org/#1}{\path{#1}}}
\providecommand{\Pubmed}[1]{\href{pmid:#1}{\path{#1}}}
\providecommand{\bibinfo}[2]{#2}
\ifx\xfnm\relax \def\xfnm[#1]{\unskip,\space#1}\fi
\bibitem[{Albiac and Kalton(2006)}]{Albiac2006}
\bibinfo{author}{Albiac, F.}, \bibinfo{author}{Kalton, N.J.},
  \bibinfo{year}{2006}.
\newblock \bibinfo{title}{Topics in {Banach} Space Theory}.
\newblock Graduate Texts in Mathematics, \bibinfo{publisher}{Springer}.
\newblock \DOIprefix\doi{10.1007/978-3-319-31557-7}.
\bibitem[{Amini et~al.(2014)Amini, Pellegrini and Rouchon}]{amini2014}
\bibinfo{author}{Amini, H.}, \bibinfo{author}{Pellegrini, C.},
  \bibinfo{author}{Rouchon, P.}, \bibinfo{year}{2014}.
\newblock \bibinfo{title}{Stability of continuous-time quantum filters with
  measurement imperfections}.
\newblock \bibinfo{journal}{Russ. J. Math. Phys.} \bibinfo{volume}{21},
  \bibinfo{pages}{297--315}.
\newblock \DOIprefix\doi{10.1134/S1061920814030029}.
\bibitem[{Amini et~al.(2013)Amini, Somaraju, Dotsenko, Sayrin, Mirrahimi and
  Rouchon}]{Amini2013}
\bibinfo{author}{Amini, H.}, \bibinfo{author}{Somaraju, R.A.},
  \bibinfo{author}{Dotsenko, I.}, \bibinfo{author}{Sayrin, C.},
  \bibinfo{author}{Mirrahimi, M.}, \bibinfo{author}{Rouchon, P.},
  \bibinfo{year}{2013}.
\newblock \bibinfo{title}{Feedback stabilization of discrete-time quantum
  systems subject to non-demolition measurements with imperfections and
  delays}.
\newblock \bibinfo{journal}{Automatica} \bibinfo{volume}{49},
  \bibinfo{pages}{2683--2692}.
\newblock \DOIprefix\doi{10.1016/j.automatica.2013.06.012}.
\bibitem[{Arnold et~al.(2008)Arnold, Fagnola and Neumann}]{Arnold2008}
\bibinfo{author}{Arnold, A.}, \bibinfo{author}{Fagnola, F.},
  \bibinfo{author}{Neumann, L.}, \bibinfo{year}{2008}.
\newblock \bibinfo{title}{Quantum {Fokker}-{Planck} models: The {Linblad} and
  {Wigner} approaches}, in: \bibinfo{booktitle}{Quantum Probability and Related
  Topics}, \bibinfo{publisher}{World Scientific}. pp. \bibinfo{pages}{23--48}.
\newblock \DOIprefix\doi{10.1142/9789812835277_0003}.
\bibitem[{Arnold et~al.(2012)Arnold, Gamba, Gualdani, Mischler, Mouhot and
  Sparber}]{ARNOLD2012}
\bibinfo{author}{Arnold, A.}, \bibinfo{author}{Gamba, I.},
  \bibinfo{author}{Gualdani, M.P.}, \bibinfo{author}{Mischler, S.},
  \bibinfo{author}{Mouhot, C.}, \bibinfo{author}{Sparber, C.},
  \bibinfo{year}{2012}.
\newblock \bibinfo{title}{The {Wigner-Fokker-Planck} equation: stationary
  states and large time behavior}.
\newblock \bibinfo{journal}{Math. Models \& Methods Appl. Sci.}
  \bibinfo{volume}{22}, \bibinfo{pages}{1250034}.
\newblock \DOIprefix\doi{10.1142/s0218202512500340}.
\bibitem[{Azouit et~al.(2015)Azouit, Sarlette and Rouchon}]{Azouit2015}
\bibinfo{author}{Azouit, R.}, \bibinfo{author}{Sarlette, A.},
  \bibinfo{author}{Rouchon, P.}, \bibinfo{year}{2015}.
\newblock \bibinfo{title}{Convergence and adiabatic elimination for a driven
  dissipative quantum harmonic oscillator}, in: \bibinfo{booktitle}{2015 54th
  {IEEE} Conference on Decision and Control ({CDC})},
  \bibinfo{publisher}{{IEEE}}.
\newblock \DOIprefix\doi{10.1109/cdc.2015.7403235}.
\bibitem[{Belavkin(1999)}]{Belavkin1999}
\bibinfo{author}{Belavkin, V.P.}, \bibinfo{year}{1999}.
\newblock \bibinfo{title}{Measurement, filtering and control in quantum open
  dynamical systems}.
\newblock \bibinfo{journal}{Rep. Math. Phys.} \bibinfo{volume}{43},
  \bibinfo{pages}{A405--A425}.
\newblock \DOIprefix\doi{10.1016/s0034-4877(00)86386-7}.
\bibitem[{Benoist et~al.(2017)Benoist, Pellegrini and Ticozzi}]{Benoist2017}
\bibinfo{author}{Benoist, T.}, \bibinfo{author}{Pellegrini, C.},
  \bibinfo{author}{Ticozzi, F.}, \bibinfo{year}{2017}.
\newblock \bibinfo{title}{Exponential stability of subspaces for quantum
  stochastic master equations}.
\newblock \bibinfo{journal}{Ann. Henri Poincar{\'{e}}}
  \DOIprefix\doi{10.1007/s00023-017-0556-3}.
\bibitem[{Bouten et~al.(2007)Bouten, Handel and James}]{bouten2007introduction}
\bibinfo{author}{Bouten, L.}, \bibinfo{author}{Handel, R.V.},
  \bibinfo{author}{James, M.R.}, \bibinfo{year}{2007}.
\newblock \bibinfo{title}{An introduction to quantum filtering}.
\newblock \bibinfo{journal}{{SIAM} J. Contr. \& Optim.} \bibinfo{volume}{46},
  \bibinfo{pages}{2199--2241}.
\newblock \DOIprefix\doi{10.1137/060651239}.
\bibitem[{Camalet(2013)}]{Camalet2013}
\bibinfo{author}{Camalet, S.}, \bibinfo{year}{2013}.
\newblock \bibinfo{title}{Steady {Schr\"{o}dinger} cat state of a driven
  {Ising} chain}.
\newblock \bibinfo{journal}{Eur. Phys. J. B} \bibinfo{volume}{86}.
\newblock \DOIprefix\doi{10.1140/epjb/e2013-31010-0}.
\bibitem[{Chang(2014)}]{Chang2014}
\bibinfo{author}{Chang, M.H.}, \bibinfo{year}{2014}.
\newblock \bibinfo{title}{A survey on invariance and ergodicity of quantum
  {Markov} semigroups}.
\newblock \bibinfo{journal}{Stoch. Anal. \& Appl.} \bibinfo{volume}{32},
  \bibinfo{pages}{480--554}.
\newblock \DOIprefix\doi{10.1080/07362994.2014.897136}.
\bibitem[{Chang(2015)}]{Chang2015}
\bibinfo{author}{Chang, M.H.}, \bibinfo{year}{2015}.
\newblock \bibinfo{title}{Quantum Stochastics}.
\newblock \bibinfo{publisher}{Cambridge University Press}.
\newblock \DOIprefix\doi{10.1017/cbo9781107706545}.
\bibitem[{Chen et~al.(2003)Chen, Gu and Kharitonov}]{Chen2003}
\bibinfo{author}{Chen, J.}, \bibinfo{author}{Gu, K.},
  \bibinfo{author}{Kharitonov, V.L.}, \bibinfo{year}{2003}.
\newblock \bibinfo{title}{Stability of Time-Delay Systems}.
\newblock \bibinfo{publisher}{Birkh\"{a}user}.
\newblock \DOIprefix\doi{10.1007/978-1-4612-0039-0}.
\bibitem[{Conway(1985)}]{Conway1985}
\bibinfo{author}{Conway, J.B.}, \bibinfo{year}{1985}.
\newblock \bibinfo{title}{A Course in Functional Analysis}.
  volume~\bibinfo{volume}{96} of \textit{\bibinfo{series}{Graduate Texts in
  Mathematics}}.
\newblock \bibinfo{publisher}{Springer New York}.
\newblock \DOIprefix\doi{10.1007/978-1-4757-3828-5}.
\bibitem[{Conway(2000)}]{conway2000course}
\bibinfo{author}{Conway, J.B.}, \bibinfo{year}{2000}.
\newblock \bibinfo{title}{A Course in Operator Theory}.
  volume~\bibinfo{volume}{2}.
\newblock \bibinfo{publisher}{American Mathematical Society}.
\newblock \DOIprefix\doi{10.1090/gsm/021}.
\bibitem[{Davies(1976)}]{davies1976quantum}
\bibinfo{author}{Davies, E.B.}, \bibinfo{year}{1976}.
\newblock \bibinfo{title}{Quantum Theory of Open Systems}.
\newblock \bibinfo{publisher}{Academic Press}.
\bibitem[{Deschamps et~al.(2016)Deschamps, Fagnola, Sasso and
  Umanit{\`{a}}}]{Deschamps2016}
\bibinfo{author}{Deschamps, J.}, \bibinfo{author}{Fagnola, F.},
  \bibinfo{author}{Sasso, E.}, \bibinfo{author}{Umanit{\`{a}}, V.},
  \bibinfo{year}{2016}.
\newblock \bibinfo{title}{Structure of uniformly continuous quantum {Markov}
  semigroups}.
\newblock \bibinfo{journal}{Rev. Math. Phys.} \bibinfo{volume}{28},
  \bibinfo{pages}{1650003}.
\newblock \DOIprefix\doi{10.1142/s0129055x16500033}.
\bibitem[{Dong and Petersen(2010)}]{Dong2010}
\bibinfo{author}{Dong, D.}, \bibinfo{author}{Petersen, I.},
  \bibinfo{year}{2010}.
\newblock \bibinfo{title}{Quantum control theory and applications: a survey}.
\newblock \bibinfo{journal}{{IET} Contr. Theor. {\&} Appl.}
  \bibinfo{volume}{4}, \bibinfo{pages}{2651--2671}.
\newblock \DOIprefix\doi{10.1049/iet-cta.2009.0508}.
\bibitem[{Emzir(2018)}]{Emzir2018}
\bibinfo{author}{Emzir, M.F.}, \bibinfo{year}{2018}.
\newblock \bibinfo{title}{Topics on filtering and coherent control of nonlinear
  quantum systems}.
\newblock Ph.D. thesis. UNSW.
\bibitem[{Emzir et~al.(2017)Emzir, Woolley and Petersen}]{Emzir2017a}
\bibinfo{author}{Emzir, M.F.}, \bibinfo{author}{Woolley, M.J.},
  \bibinfo{author}{Petersen, I.R.}, \bibinfo{year}{2017}.
\newblock \bibinfo{title}{{Lyapunov} stability analysis for invariant states of
  quantum systems}, in: \bibinfo{booktitle}{Proc of 2017 {IEEE} 56th Conf. on
  Dec. and Contr. ({CDC})}.
\newblock \DOIprefix\doi{10.1109/cdc.2017.8264475}.
\bibitem[{Emzir et~al.(2018)Emzir, Woolley and Petersen}]{Emzir2018a}
\bibinfo{author}{Emzir, M.F.}, \bibinfo{author}{Woolley, M.J.},
  \bibinfo{author}{Petersen, I.R.}, \bibinfo{year}{2018}.
\newblock \bibinfo{title}{On physical realizability of nonlinear quantum
  stochastic differential equations}.
\newblock \bibinfo{journal}{Automatica} \bibinfo{volume}{95},
  \bibinfo{pages}{254--265}.
\newblock \DOIprefix\doi{10.1016/j.automatica.2018.05.011}.
\bibitem[{Fagnola(1999)}]{fagnola1999quantum}
\bibinfo{author}{Fagnola, F.}, \bibinfo{year}{1999}.
\newblock \bibinfo{title}{Quantum {Markov} semigroups and quantum flows}.
\newblock \bibinfo{journal}{Proyecciones} \bibinfo{volume}{18},
  \bibinfo{pages}{1--144}.
\bibitem[{Fagnola and Rebolledo(2001)}]{Fagnola2001}
\bibinfo{author}{Fagnola, F.}, \bibinfo{author}{Rebolledo, R.},
  \bibinfo{year}{2001}.
\newblock \bibinfo{title}{On the existence of stationary states for quantum
  dynamical semigroups}.
\newblock \bibinfo{journal}{J. Math. Phys.} \bibinfo{volume}{42},
  \bibinfo{pages}{1296}.
\newblock \DOIprefix\doi{10.1063/1.1340870}.
\bibitem[{Fagnola and Rebolledo(2002)}]{Fagnola2002}
\bibinfo{author}{Fagnola, F.}, \bibinfo{author}{Rebolledo, R.},
  \bibinfo{year}{2002}.
\newblock \bibinfo{title}{Lectures on the qualitative analysis of quantum
  {Markov} semigroups}, in: \bibinfo{booktitle}{Quantum Interacting Particle
  Systems}. \bibinfo{publisher}{World Scientific}, pp.
  \bibinfo{pages}{197--239}.
\newblock \DOIprefix\doi{10.1142/9789812776853_0002}.
\bibitem[{Fagnola and Rebolledo(2003)}]{Fagnola2003}
\bibinfo{author}{Fagnola, F.}, \bibinfo{author}{Rebolledo, R.},
  \bibinfo{year}{2003}.
\newblock \bibinfo{title}{Quantum {Markov} semigroups and their stationary
  states}, in: \bibinfo{booktitle}{Stochastic Analysis and Mathematical Physics
  {II}}. \bibinfo{publisher}{Springer Science Business Media}, pp.
  \bibinfo{pages}{77--128}.
\newblock \DOIprefix\doi{10.1007/978-3-0348-8018-3_6}.
\bibitem[{Fagnola and Rebolledo(2006)}]{Fagnola2006}
\bibinfo{author}{Fagnola, F.}, \bibinfo{author}{Rebolledo, R.},
  \bibinfo{year}{2006}.
\newblock \bibinfo{title}{Notes on the qualitative behaviour of quantum
  {Markov} semigroups}, in: \bibinfo{booktitle}{Open Quantum Systems {III}}.
  \bibinfo{publisher}{Springer}, pp. \bibinfo{pages}{161--205}.
\newblock \DOIprefix\doi{10.1007/3-540-33967-1_4}.
\bibitem[{Fuchs and van~de Graaf(1999)}]{Fuchs1999}
\bibinfo{author}{Fuchs, C.A.}, \bibinfo{author}{van~de Graaf, J.},
  \bibinfo{year}{1999}.
\newblock \bibinfo{title}{Cryptographic distinguishability measures for
  quantum-mechanical states}.
\newblock \bibinfo{journal}{{IEEE} Trans. Inform. Theory} \bibinfo{volume}{45},
  \bibinfo{pages}{1216--1227}.
\newblock \DOIprefix\doi{10.1109/18.761271}.
\bibitem[{Gantmacher(1959)}]{GantMacher1959}
\bibinfo{author}{Gantmacher, F.R.}, \bibinfo{year}{1959}.
\newblock \bibinfo{title}{The Theory of Matrices}. volume~\bibinfo{volume}{1}.
\newblock \bibinfo{publisher}{Chelsea Publishing Company}.
\bibitem[{Gough and James(2009)}]{gough2009series}
\bibinfo{author}{Gough, J.}, \bibinfo{author}{James, M.}, \bibinfo{year}{2009}.
\newblock \bibinfo{title}{The series product and its application to quantum
  feedforward and feedback networks}.
\newblock \bibinfo{journal}{{IEEE} Trans. Automat. Contr.}
  \bibinfo{volume}{54}, \bibinfo{pages}{2530--2544}.
\newblock \DOIprefix\doi{10.1109/tac.2009.2031205}.
\bibitem[{Hamerly and Mabuchi(2012)}]{Hamerly2012}
\bibinfo{author}{Hamerly, R.}, \bibinfo{author}{Mabuchi, H.},
  \bibinfo{year}{2012}.
\newblock \bibinfo{title}{Advantages of coherent feedback for cooling quantum
  oscillators}.
\newblock \bibinfo{journal}{Phys. Rev. Lett.} \bibinfo{volume}{109}.
\newblock \DOIprefix\doi{10.1103/physrevlett.109.173602}.
\bibitem[{Ingemar~Bengtsson(2017)}]{IngemarBengtsson2017}
\bibinfo{author}{Ingemar~Bengtsson, K.Z.}, \bibinfo{year}{2017}.
\newblock \bibinfo{title}{Geometry of Quantum States}.
\newblock \bibinfo{publisher}{Cambridge University Pr.}
\bibitem[{James et~al.(2008)James, Nurdin and Petersen}]{James2008}
\bibinfo{author}{James, M.}, \bibinfo{author}{Nurdin, H.},
  \bibinfo{author}{Petersen, I.}, \bibinfo{year}{2008}.
\newblock \bibinfo{title}{{$\mathcal{H}_{\infty}$} control of linear quantum
  stochastic systems}.
\newblock \bibinfo{journal}{{IEEE} Trans. Automat. Contr.}
  \bibinfo{volume}{53}, \bibinfo{pages}{1787--1803}.
\newblock \DOIprefix\doi{10.1109/tac.2008.929378}.
\bibitem[{James and Gough(2010)}]{James2010}
\bibinfo{author}{James, M.R.}, \bibinfo{author}{Gough, J.E.},
  \bibinfo{year}{2010}.
\newblock \bibinfo{title}{Quantum dissipative systems and feedback control
  design by interconnection}.
\newblock \bibinfo{journal}{{IEEE} Trans. Automat. Contr.}
  \bibinfo{volume}{55}, \bibinfo{pages}{1806--1821}.
\newblock \DOIprefix\doi{10.1109/tac.2010.2046067}.
\bibitem[{Jungers et~al.(2017)Jungers, Ahmadi, Parrilo and
  Roozbehani}]{Jungers2017}
\bibinfo{author}{Jungers, R.}, \bibinfo{author}{Ahmadi, A.A.},
  \bibinfo{author}{Parrilo, P.A.}, \bibinfo{author}{Roozbehani, M.},
  \bibinfo{year}{2017}.
\newblock \bibinfo{title}{A characterization of {Lyapunov} inequalities for
  stability of switched systems}.
\newblock \bibinfo{journal}{{IEEE} Trans. Automat. Contr.} ,
  \bibinfo{pages}{1--1}\DOIprefix\doi{10.1109/tac.2017.2671345}.
\bibitem[{Khalil(2002)}]{khalil2002nonlinear}
\bibinfo{author}{Khalil, H.}, \bibinfo{year}{2002}.
\newblock \bibinfo{title}{Nonlinear Systems}.
\newblock \bibinfo{publisher}{Prentice Hall}.
\bibitem[{Kharitonov(1999)}]{Kharitonov1999}
\bibinfo{author}{Kharitonov, V.L.}, \bibinfo{year}{1999}.
\newblock \bibinfo{title}{Robust stability analysis of time delay systems: A
  survey}.
\newblock \bibinfo{journal}{Annual Reviews in Control} \bibinfo{volume}{23},
  \bibinfo{pages}{185--196}.
\newblock \DOIprefix\doi{10.1016/s1367-5788(99)90087-1}.
\bibitem[{Khasminskii(2011)}]{khasminskii2011stochastic}
\bibinfo{author}{Khasminskii, R.}, \bibinfo{year}{2011}.
\newblock \bibinfo{title}{Stochastic Stability of Differential Equations}.
  volume~\bibinfo{volume}{66}.
\newblock \bibinfo{publisher}{Springer}.
\newblock \DOIprefix\doi{10.1007/978-3-642-23280-0}.
\bibitem[{Kozin(1969)}]{Kozin1969}
\bibinfo{author}{Kozin, F.}, \bibinfo{year}{1969}.
\newblock \bibinfo{title}{A survey of stability of stochastic systems}.
\newblock \bibinfo{journal}{Automatica} \bibinfo{volume}{5},
  \bibinfo{pages}{95--112}.
\newblock \DOIprefix\doi{10.1016/0005-1098(69)90060-0}.
\bibitem[{Krasovsky(1963)}]{Krasovsky1963}
\bibinfo{author}{Krasovsky, N.N.}, \bibinfo{year}{1963}.
\newblock \bibinfo{title}{Stability of Motion: Applications of {Lyapunov}'s
  Second Method to Differential Systems and Equations with Delay}.
\newblock Studies in mathematical analysis and related topics,
  \bibinfo{publisher}{Stanford University Press}.
\bibitem[{Leghtas et~al.(2015)Leghtas, Touzard, Pop, Kou, Vlastakis, Petrenko,
  Sliwa, Narla, Shankar, Hatridge, Reagor, Frunzio, Schoelkopf, Mirrahimi and
  Devoret}]{Leghtas2015}
\bibinfo{author}{Leghtas, Z.}, \bibinfo{author}{Touzard, S.},
  \bibinfo{author}{Pop, I.M.}, \bibinfo{author}{Kou, A.},
  \bibinfo{author}{Vlastakis, B.}, \bibinfo{author}{Petrenko, A.},
  \bibinfo{author}{Sliwa, K.M.}, \bibinfo{author}{Narla, A.},
  \bibinfo{author}{Shankar, S.}, \bibinfo{author}{Hatridge, M.J.},
  \bibinfo{author}{Reagor, M.}, \bibinfo{author}{Frunzio, L.},
  \bibinfo{author}{Schoelkopf, R.J.}, \bibinfo{author}{Mirrahimi, M.},
  \bibinfo{author}{Devoret, M.H.}, \bibinfo{year}{2015}.
\newblock \bibinfo{title}{Confining the state of light to a quantum manifold by
  engineered two-photon loss}.
\newblock \bibinfo{journal}{Science} \bibinfo{volume}{347},
  \bibinfo{pages}{853--857}.
\newblock \DOIprefix\doi{10.1126/science.aaa2085}.
\bibitem[{Leghtas et~al.(2013)Leghtas, Vool, Shankar, Hatridge, Girvin, Devoret
  and Mirrahimi}]{Leghtas2013}
\bibinfo{author}{Leghtas, Z.}, \bibinfo{author}{Vool, U.},
  \bibinfo{author}{Shankar, S.}, \bibinfo{author}{Hatridge, M.},
  \bibinfo{author}{Girvin, S.M.}, \bibinfo{author}{Devoret, M.H.},
  \bibinfo{author}{Mirrahimi, M.}, \bibinfo{year}{2013}.
\newblock \bibinfo{title}{Stabilizing a {Bell} state of two superconducting
  qubits by dissipation engineering}.
\newblock \bibinfo{journal}{Phys. Rev. A: At., Mol., Opt. Phys.}
  \bibinfo{volume}{88}.
\newblock \DOIprefix\doi{10.1103/physreva.88.023849}.
\bibitem[{Lin and Antsaklis(2009)}]{Lin2009}
\bibinfo{author}{Lin, H.}, \bibinfo{author}{Antsaklis, P.J.},
  \bibinfo{year}{2009}.
\newblock \bibinfo{title}{Stability and stabilizability of switched linear
  systems: A survey of recent results}.
\newblock \bibinfo{journal}{{IEEE} Trans. Automat. Contr.}
  \bibinfo{volume}{54}, \bibinfo{pages}{308--322}.
\newblock \DOIprefix\doi{10.1109/tac.2008.2012009}.
\bibitem[{Liu and Sun(2016)}]{Liu2016}
\bibinfo{author}{Liu, K.Z.}, \bibinfo{author}{Sun, X.M.}, \bibinfo{year}{2016}.
\newblock \bibinfo{title}{{Razumikhin}-type theorems for hybrid system with
  memory}.
\newblock \bibinfo{journal}{Automatica} \bibinfo{volume}{71},
  \bibinfo{pages}{72--77}.
\newblock \DOIprefix\doi{10.1016/j.automatica.2016.04.038}.
\bibitem[{Lyapunov(1892)}]{Lyapunov1892}
\bibinfo{author}{Lyapunov, A.M.}, \bibinfo{year}{1892}.
\newblock \bibinfo{title}{The general problem of motion stability}.
\newblock \bibinfo{journal}{Ann. Math. Stud.} \bibinfo{volume}{17}.
\bibitem[{Maalouf and Petersen(2011)}]{Maalouf2011}
\bibinfo{author}{Maalouf, A.I.}, \bibinfo{author}{Petersen, I.R.},
  \bibinfo{year}{2011}.
\newblock \bibinfo{title}{Coherent {$\mathcal{H}_{\infty}$} control for a class
  of annihilation operator linear quantum systems}.
\newblock \bibinfo{journal}{{IEEE} Trans. Automat. Contr.}
  \bibinfo{volume}{56}, \bibinfo{pages}{309--319}.
\newblock \DOIprefix\doi{10.1109/tac.2010.2052942}.
\bibitem[{Mamaev et~al.(2018)Mamaev, Govia and Clerk}]{Mamaev2017}
\bibinfo{author}{Mamaev, M.}, \bibinfo{author}{Govia, L.C.G.},
  \bibinfo{author}{Clerk, A.A.}, \bibinfo{year}{2018}.
\newblock \bibinfo{title}{Dissipative stabilization of entangled cat states
  using a driven {Bose-Hubbard} dimer}.
\newblock \bibinfo{journal}{Quantum} \bibinfo{volume}{2}, \bibinfo{pages}{58}.
\bibitem[{Mazenc and Malisoff(2017)}]{Mazenc2017}
\bibinfo{author}{Mazenc, F.}, \bibinfo{author}{Malisoff, M.},
  \bibinfo{year}{2017}.
\newblock \bibinfo{title}{Extensions of {Razumikhin's} theorem and
  {Lyapunov-Krasovskii} functional constructions for time-varying systems with
  delay}.
\newblock \bibinfo{journal}{Automatica} \bibinfo{volume}{78},
  \bibinfo{pages}{1--13}.
\newblock \DOIprefix\doi{10.1016/j.automatica.2016.12.005}.
\bibitem[{Megginson(1998)}]{Megginson1998}
\bibinfo{author}{Megginson, R.E.}, \bibinfo{year}{1998}.
\newblock \bibinfo{title}{An Introduction to {Banach} Space Theory}.
\newblock \bibinfo{publisher}{Springer}.
\bibitem[{Mendelson(2003)}]{Mendelson2003}
\bibinfo{author}{Mendelson, B.}, \bibinfo{year}{2003}.
\newblock \bibinfo{title}{Introduction to Topology}.
\newblock \bibinfo{publisher}{Dover Publications Inc.}
\bibitem[{Meyer(1993)}]{Meyer1993}
\bibinfo{author}{Meyer, P.A.}, \bibinfo{year}{1993}.
\newblock \bibinfo{title}{Quantum Probability for Probabilists}.
\newblock \bibinfo{publisher}{Springer Berlin Heidelberg}.
\newblock \DOIprefix\doi{10.1007/978-3-662-21558-6}.
\bibitem[{Michael A.~Nielsen(2001)}]{MichaelA.Nielsen2001}
\bibinfo{author}{Michael A.~Nielsen, I.L.C.}, \bibinfo{year}{2001}.
\newblock \bibinfo{title}{Quantum Computation and Quantum Information}.
\newblock \bibinfo{publisher}{Cambridge University Pr.}
\bibitem[{Mirrahimi et~al.(2014)Mirrahimi, Leghtas, Albert, Touzard,
  Schoelkopf, Jiang and Devoret}]{Mirrahimi2014}
\bibinfo{author}{Mirrahimi, M.}, \bibinfo{author}{Leghtas, Z.},
  \bibinfo{author}{Albert, V.V.}, \bibinfo{author}{Touzard, S.},
  \bibinfo{author}{Schoelkopf, R.J.}, \bibinfo{author}{Jiang, L.},
  \bibinfo{author}{Devoret, M.H.}, \bibinfo{year}{2014}.
\newblock \bibinfo{title}{Dynamically protected cat-qubits: a new paradigm for
  universal quantum computation}.
\newblock \bibinfo{journal}{New J. Phys.} \bibinfo{volume}{16},
  \bibinfo{pages}{045014}.
\newblock \DOIprefix\doi{10.1088/1367-2630/16/4/045014}.
\bibitem[{Morris(2016)}]{Morris2016}
\bibinfo{author}{Morris, S.A.}, \bibinfo{year}{2016}.
\newblock \bibinfo{title}{Topology Without Tears}.
\newblock \URLprefix \url{http://www.topologywithouttears.net/topbook.pdf}.
\bibitem[{Nigro(2018)}]{Nigro2019}
\bibinfo{author}{Nigro, D.}, \bibinfo{year}{2018}.
\newblock \bibinfo{title}{On the uniqueness of the steady-state solution of the
  lindblad-gorini-kossakowski-sudarshan equation}.
\newblock \bibinfo{journal}{J. Stat. Mech: Theor \& Expr.}
  \bibinfo{volume}{2019}, \bibinfo{pages}{043202}.
\newblock \DOIprefix\doi{10.1088/1742-5468/ab0c1c}.
\bibitem[{Nurdin et~al.(2009)Nurdin, James and Petersen}]{Nurdin2009}
\bibinfo{author}{Nurdin, H.I.}, \bibinfo{author}{James, M.R.},
  \bibinfo{author}{Petersen, I.R.}, \bibinfo{year}{2009}.
\newblock \bibinfo{title}{Coherent quantum {LQG} control}.
\newblock \bibinfo{journal}{Automatica} \bibinfo{volume}{45},
  \bibinfo{pages}{1837--1846}.
\newblock \DOIprefix\doi{10.1016/j.automatica.2009.04.018}.
\bibitem[{Pan et~al.(2014)Pan, Amini, Miao, Gough, Ugrinovskii and
  James}]{Pan2014}
\bibinfo{author}{Pan, Y.}, \bibinfo{author}{Amini, H.}, \bibinfo{author}{Miao,
  Z.}, \bibinfo{author}{Gough, J.}, \bibinfo{author}{Ugrinovskii, V.},
  \bibinfo{author}{James, M.R.}, \bibinfo{year}{2014}.
\newblock \bibinfo{title}{{Heisenberg} picture approach to the stability of
  quantum {Markov} systems}.
\newblock \bibinfo{journal}{J. Math. Phys.} \bibinfo{volume}{55},
  \bibinfo{pages}{06270101 -- 06270116}.
\newblock \DOIprefix\doi{http://dx.doi.org/10.1063/1.4884300}.
\bibitem[{Pan et~al.(2016)Pan, Ugrinovskii and James}]{Pan2016}
\bibinfo{author}{Pan, Y.}, \bibinfo{author}{Ugrinovskii, V.},
  \bibinfo{author}{James, M.R.}, \bibinfo{year}{2016}.
\newblock \bibinfo{title}{Ground-state stabilization of quantum finite-level
  systems by dissipation}.
\newblock \bibinfo{journal}{Automatica} \bibinfo{volume}{65},
  \bibinfo{pages}{147 -- 159}.
\newblock \DOIprefix\doi{10.1016/j.automatica.2015.11.041}.
\bibitem[{Parthasarathy(1967)}]{Parthasarathy1967}
\bibinfo{author}{Parthasarathy, K.R.}, \bibinfo{year}{1967}.
\newblock \bibinfo{title}{Probability Measures on Metric Spaces}.
\newblock A Series of Monographs and Textbooks, \bibinfo{publisher}{Academic
  Press}.
\newblock \DOIprefix\doi{10.1016/B978-1-4832-0022-4.50001-6}.
\bibitem[{Parthasarathy(1992)}]{Parthasarathy1992}
\bibinfo{author}{Parthasarathy, K.R.}, \bibinfo{year}{1992}.
\newblock \bibinfo{title}{An Introduction to Quantum Stochastic Calculus}.
\newblock Modern Birkh{\"a}user classics, \bibinfo{publisher}{Springer}.
\newblock \DOIprefix\doi{10.1007/978-3-0348-0566-7}.
\bibitem[{Prokhorov(1956)}]{Prokhorov1956}
\bibinfo{author}{Prokhorov, Y.V.}, \bibinfo{year}{1956}.
\newblock \bibinfo{title}{Convergence of random processes and limit theorems in
  probability theory}.
\newblock \bibinfo{journal}{Theor. of Prob. \& Appl.} \bibinfo{volume}{1},
  \bibinfo{pages}{157--214}.
\newblock \DOIprefix\doi{10.1137/1101016}.
\bibitem[{Reed and Simon(1972)}]{Reed1972}
\bibinfo{author}{Reed, M.}, \bibinfo{author}{Simon, B.}, \bibinfo{year}{1972}.
\newblock \bibinfo{title}{Methods of Modern Mathematical Physics}.
  volume~\bibinfo{volume}{1}.
\newblock \bibinfo{publisher}{Academic press}.
\bibitem[{Salle(1960)}]{Salle1960}
\bibinfo{author}{Salle, J.L.}, \bibinfo{year}{1960}.
\newblock \bibinfo{title}{Some extensions of {Liapunov's} second method}.
\newblock \bibinfo{journal}{{IRE} Transactions on Circuit Theory}
  \bibinfo{volume}{7}, \bibinfo{pages}{520--527}.
\newblock \DOIprefix\doi{10.1109/tct.1960.1086720}.
\bibitem[{Salle and Lefschetz(1961)}]{Salle1961}
\bibinfo{author}{Salle, J.L.}, \bibinfo{author}{Lefschetz, S.},
  \bibinfo{year}{1961}.
\newblock \bibinfo{title}{Stability by {Liapunov}'s Direct Method: With
  Applications}.
\newblock Mathematics in Science and Engineering, \bibinfo{publisher}{Academic
  Press}.
\bibitem[{Schirmer and Wang(2010)}]{Schirmer2010}
\bibinfo{author}{Schirmer, S.G.}, \bibinfo{author}{Wang, X.},
  \bibinfo{year}{2010}.
\newblock \bibinfo{title}{Stabilizing open quantum systems by {Markovian}
  reservoir engineering}.
\newblock \bibinfo{journal}{Phys. Rev. A: At., Mol., Opt. Phys.}
  \bibinfo{volume}{81}.
\newblock \DOIprefix\doi{10.1103/physreva.81.062306}.
\bibitem[{Somaraju et~al.(2013)Somaraju, Mirrahimi and Rouchon}]{somaraju2013}
\bibinfo{author}{Somaraju, R.}, \bibinfo{author}{Mirrahimi, M.},
  \bibinfo{author}{Rouchon, P.}, \bibinfo{year}{2013}.
\newblock \bibinfo{title}{Approximate stabilization of an infinite dimensional
  quantum stochastic system}.
\newblock \bibinfo{journal}{Rev. Math. Phys.} \bibinfo{volume}{25},
  \bibinfo{pages}{1350001}.
\newblock \DOIprefix\doi{10.1142/s0129055x13500013}.
\bibitem[{Sutherland(2009)}]{Sutherland2009}
\bibinfo{author}{Sutherland, W.A.}, \bibinfo{year}{2009}.
\newblock \bibinfo{title}{Introduction to Metric and Topological Spaces}.
\newblock \bibinfo{publisher}{Oxford University Press}.
\bibitem[{Ticozzi and Viola(2012)}]{Ticozzi2012}
\bibinfo{author}{Ticozzi, F.}, \bibinfo{author}{Viola, L.},
  \bibinfo{year}{2012}.
\newblock \bibinfo{title}{Stabilizing entangled states with quasi-local quantum
  dynamical semigroups}.
\newblock \bibinfo{journal}{Philos. Trans. R. Soc. London, Ser. A}
  \bibinfo{volume}{370}, \bibinfo{pages}{5259--5269}.
\newblock \DOIprefix\doi{10.1098/rsta.2011.0485}.
\bibitem[{Wiseman and Milburn(1993)}]{Wiseman1993}
\bibinfo{author}{Wiseman, H.M.}, \bibinfo{author}{Milburn, G.J.},
  \bibinfo{year}{1993}.
\newblock \bibinfo{title}{Quantum theory of optical feedback via homodyne
  detection}.
\newblock \bibinfo{journal}{Phys. Rev. Lett.} \bibinfo{volume}{70},
  \bibinfo{pages}{548--551}.
\newblock \DOIprefix\doi{10.1103/physrevlett.70.548}.
\bibitem[{Wiseman and Milburn(1994)}]{Wiseman1994}
\bibinfo{author}{Wiseman, H.M.}, \bibinfo{author}{Milburn, G.J.},
  \bibinfo{year}{1994}.
\newblock \bibinfo{title}{All-optical versus electro-optical quantum-limited
  feedback}.
\newblock \bibinfo{journal}{Phys. Rev. A: At., Mol., Opt. Phys.}
  \bibinfo{volume}{49}, \bibinfo{pages}{4110--4125}.
\newblock \DOIprefix\doi{10.1103/PhysRevA.49.4110}.
\bibitem[{Wiseman and Milburn(2010)}]{wiseman2010quantum}
\bibinfo{author}{Wiseman, H.M.}, \bibinfo{author}{Milburn, G.J.},
  \bibinfo{year}{2010}.
\newblock \bibinfo{title}{Quantum Measurement and Control}.
\newblock \bibinfo{publisher}{Cambridge University Press}.
\newblock \DOIprefix\doi{10.1017/cbo9780511813948.002}.
\bibitem[{Wojtaszczyk(1991)}]{Wojtaszczyk1991}
\bibinfo{author}{Wojtaszczyk, P.}, \bibinfo{year}{1991}.
\newblock \bibinfo{title}{{Banach} Spaces For Analysts}.
\newblock \bibinfo{publisher}{Cambridge University Press}.
\newblock \DOIprefix\doi{10.1017/cbo9780511608735}.
\bibitem[{Yamamoto(2012)}]{Yamamoto2012}
\bibinfo{author}{Yamamoto, N.}, \bibinfo{year}{2012}.
\newblock \bibinfo{title}{Pure {Gaussian} state generation via dissipation: a
  quantum stochastic differential equation approach}.
\newblock \bibinfo{journal}{Philos. Trans. R. Soc. London, Ser. A}
  \bibinfo{volume}{370}, \bibinfo{pages}{5324--5337}.
\newblock \DOIprefix\doi{10.1098/rsta.2011.0529}.
\bibitem[{Yamamoto(2014)}]{Yamamoto2014}
\bibinfo{author}{Yamamoto, N.}, \bibinfo{year}{2014}.
\newblock \bibinfo{title}{Coherent versus measurement feedback: Linear systems
  theory for quantum information}.
\newblock \bibinfo{journal}{Phys. Rev. X} \bibinfo{volume}{4}.
\newblock \DOIprefix\doi{10.1103/physrevx.4.041029}.
\bibitem[{Zhang et~al.(2012)Zhang, Wu, xi~Liu, Li and Tarn}]{Zhang2012}
\bibinfo{author}{Zhang, J.}, \bibinfo{author}{Wu, R.B.},
  \bibinfo{author}{xi~Liu, Y.}, \bibinfo{author}{Li, C.W.},
  \bibinfo{author}{Tarn, T.J.}, \bibinfo{year}{2012}.
\newblock \bibinfo{title}{Quantum coherent nonlinear feedback with applications
  to quantum optics on chip}.
\newblock \bibinfo{journal}{{IEEE} Trans. Automat. Contr.}
  \bibinfo{volume}{57}, \bibinfo{pages}{1997--2008}.
\newblock \DOIprefix\doi{10.1109/tac.2012.2195871}.
\bibitem[{Zhou and Egorov(2016)}]{Zhou2016}
\bibinfo{author}{Zhou, B.}, \bibinfo{author}{Egorov, A.V.},
  \bibinfo{year}{2016}.
\newblock \bibinfo{title}{{Razumikhin} and {Krasovskii} stability theorems for
  time-varying time-delay systems}.
\newblock \bibinfo{journal}{Automatica} \bibinfo{volume}{71},
  \bibinfo{pages}{281--291}.
\newblock \DOIprefix\doi{10.1016/j.automatica.2016.04.048}.

\end{thebibliography}

\end{document}